\documentclass[aos,nameyear,dvips]{arximspdf}
\usepackage{mathbh}
\usepackage{graphicx}

\doi{10.1214/009053607000000820}
\volume{36}
\issue{3}
\pubyear{2008}
\firstpage{1464}
\lastpage{1507}

\makeatletter
\newtheorem{defn}{Definition}[section]
\newtheorem{lem}[defn]{Lemma}
\newtheorem{lemm}{Lemma}
\newtheorem{lemma}{Lemma}
\newtheorem{lemmas}{Lemma}
\newtheorem{theorem}[defn]{Theorem}
\newtheorem{theo}{Theorem}
\newtheorem{theor}{Theorem}
\newtheorem{cor}[defn]{Corollary}
\newtheorem{corol}{Corollary}
\newproclaim{assu}[defn]{Assumption}
\newproclaim{assumpt}{Assumption}
\newproclaim{exa}[defn]{Example}
\newproclaim{rest}[defn]{Restriction}
\let \epsilon=\varepsilon

\makeatother

\begin{document}
\begin{frontmatter}

\title{Statistical modeling of causal effects in continuous
  time\protect\thanksref{T1}}
  \thankstext{T1}{Supported in part by NIH Grants R37AI032475 and R01AI51164.}
\runtitle{Causal effects in continuous time}

\begin{aug}
\author[A]{\fnms{Judith J.} \snm{Lok}\corref{}\ead[label=e1]{jlok@hsph.harvard.edu}}
\runauthor{J. J. Lok}
\affiliation{Free University of Amsterdam
and\\
Harvard School of Public Health}
\address[A]{Department of Biostatistics\\
Harvard School of Public Health\\
Building II, Room 409\\
655 Huntington Avenue\\
Boston, Massachusetts 02115\\
USA\\
\printead{e1}} 
\end{aug}

\received{\smonth{6} \syear{2004}}
\revised{\smonth{3} \syear{2007}}

\begin{abstract}
This article studies the estimation of the causal effect of a
  time-varying treatment on time-to-an-event or on some other
  continuously distributed outcome. The paper applies to the situation where
  treatment is repeatedly adapted to time-dependent patient
  characteristics. The treatment effect cannot be estimated
  by simply conditioning on these time-dependent patient
  characteristics, as they may
  themselves be indications of the treatment effect. This
  time-dependent confounding is common in observational
  studies. Robins [(1992) \textit{Biometrika} \textbf{79} 321--334,
  (1998b) \textit{Encyclopedia of Biostatistics} \textbf{6} 4372--4389]
  has proposed the so-called structural
  nested models to estimate treatment effects in the presence of
  time-dependent confounding. In this article we provide a conceptual
  framework and formalization for structural nested models in
  continuous time. We show that the resulting estimators are consistent
  and asymptotically normal. Moreover, as conjectured in Robins [(1998b) \textit{Encyclopedia of Biostatistics}
  \textbf{6} 4372--4389],
  a test for whether treatment affects the outcome of interest can be
  performed without specifying a model for treatment effect. We
  illustrate the ideas in this article with an example.
\end{abstract}

\begin{keyword}[class=AMS]
\kwd[Primary ]{62P10}
\kwd[; secondary ]{62M99}.
\end{keyword}
\begin{keyword}
\kwd{Causality in continuous time}
\kwd{counterfactuals}
\kwd{longitudinal data}
\kwd{observational studies}.
\end{keyword}

\end{frontmatter}

\section{Introduction}
\label{intro}

Causality is a topic which nowadays receives much attention.
Statisticians, epidemiologists, biostatisticians, social scientists,
computer \mbox{scientists} [especially those in artificial intelligence, see,
e.g., \citeauthor{Pea} (\citeyear{Pea})], econometricians and
philosophers are investigating questions like ``what would have
happened if'' and ``what would happen if.'' This article discusses
estimating the effect of a time-varying treatment. As a recurring
example, this article focuses on the effect of a medical treatment
which is adapted to a patient's state during the course of time.

Large observational studies have become widely used in medical
research when data from randomized experiments are not
available. Randomized clinical trials are often expensive,
impractical, and sometimes unfeasible for ethical reasons because
treatment is withheld from some patients regardless of medical
considerations. Also, in some instances, exploratory investigations
using nonexperimental data are used before conducting a randomized
trial. In observational\vadjust{\goodbreak} studies there is no pre-specified treatment
protocol. Data are collected on patient characteristics and treatments
in the course of the normal interaction between patients and
doctors. Obviously, it is considerably more difficult to draw correct
causal conclusions from observational data than from a randomized
experiment. The main reason is the so-called confounding by indication
or selection bias. For example, doctors may prescribe more medication to
patients who are relatively unhealthy. Thus, association between
medication dose and health outcomes may arise not only from the
treatment effect but also from the way the treatment was assigned.

If this confounding by indication only takes place at the start of the
treatment, one can condition on initial patient characteristics or
covariates, such as blood pressure or number of white blood cells, in
order to remove the effect of the confounding, and get meaningful
estimates of the treatment effect. Linear regression, logistic
regression or Cox regression can be used for this purpose. However,
estimating treatment effects is more difficult if treatment decisions
after the start of the treatment are adapted to the state of the
patients in subsequent periods. Treatment might be influenced by a
patient's state in the past, which was influenced by treatment
decisions before; thus, simply conditioning on a patient's state in
the past means disregarding information on the effect of past
treatment. In such a case, even the well-known time-dependent Cox model
does not answer the question of whether, or how, treatment affects the
outcome of interest. The time-dependent Cox model studies the rate at
which some event of interest happens (e.g., the patient dying), given
past treatment- and covariate history. However, under time-dependent
confounding, past covariate values may have been influenced by
previous treatment. The net effect of treatment can thus not be
derived from just this rate; see, for example, \citet{Enc}, \citet{Kei} or
\citet{Lok}.

Structural nested models, proposed in \citet{R89}, Lok,
Gill,  van der Vaart and Robins (\citeyear{SNart}) and
Robins (\citeyear{R92,Enc}) to solve practical problems in epidemiology and
biostatistics, effectively overcome these difficulties and estimate
the effect of time-varying treatments. The main assumption underlying
these models is that all information the doctors used to make
treatment decisions, and which is predictive of the patient's
prognosis with respect to the final outcome, is available for
analysis. This assumption of ``no unmeasured confounding'' makes it
possible to distinguish between treatment effect and selection bias.
What data have to be collected to satisfy this assumption of no
unmeasured confounding is for subject matter experts to decide. All of
the past treatment- and covariate information which both
(i)~influences a doctor's treatment decisions and (ii)~is relevant for
a patient's prognosis with respect to the outcome of interest, has to
be recorded. In Section~5 of \citet{Enc} and in Section~8.1 of
Robins, Rotnitzky and Scharfstein (\citeyear{RRS}), a sensitivity analysis methodology for estimation of
structural nested models is developed that does not assume no
unmeasured confounders. Beyond treatment and covariates, the data
requirements also include the measure of an outcome of interest; for example,
 survival time, time to clinical AIDS or CD4 count after the treatment
period.

Lok et al. (\citeyear{SNart}) study structural nested models in discrete time. These
models assume that changes in the values of the covariates and
treatment decisions take place at finitely many deterministic times,
which are the same for all patients and known in advance.
Lok et al. (\citeyear{SNart}) also assume that covariates and treatment take values in
a discrete space. They indicate why it is reasonable to expect
consistency and asymptotic normality in discrete time, and they refer
to \citet{Lok} for the proofs. \citet{RJ} generalize Lok et al.
(\citeyear{SNart})
to covariates and treatment taking values in ${\mathbb R}^k$.

In this article we consider structural nested models in continuous
time, proposed in Robins (\citeyear{R92,Enc}). Structural nested models in
continuous time allow for both changes in the values of the covariates
and treatment decisions to take place at arbitrary times for different
patients. As noted in \citet{corr}, structural nested models in
continuous time assume that a short duration of treatment has only a
small effect on the distribution of the outcome of interest. The
effect of the treatment on an individual patient may be large, but
then the probability of such effect has to be small for any particular
short duration of treatment (see page~7, bottom).

This article provides a conceptual framework and mathematical
formalization of these practical methods, solving important
outstanding problems and contributing to the causality discussion,
especially for the time ordered and continuous time case.  In
particular, this article proves the conjectures in \citet{Enc} that
structural nested models in continuous time lead to estimators which
are both consistent and asymptotically normal. The proof simplifies
considerably for structural nested models in discrete time (see our
Discussion, Section~\ref{discu}). This article also proves that a test
related to the score test can be used to investigate whether treatment
affects the outcome of interest without specifying a model for the
treatment effect.

\section{Setting and notation}
\label{sets}

The setting to which structural nested models in continuous time apply
is as follows. The outcome of interest, from now on called $Y$, is a
continuous real variable. For example, the survival time of a patient,
time to clinical AIDS, or CD4 count after the treatment period. We
wish to estimate the effect of treatment on the outcome $Y$. There is
some fixed time interval $[0,\tau]$, with $\tau$ a finite
time, during which treatment and patient characteristics are observed
for each patient. We suppose that after time $\tau$ treatment is
stopped or switched to some kind of baseline treatment. In this
article we assume that there is no censoring, and that the outcome $Y$
is observed for every patient in the study. See, for example, \citet{Enc},
Hern\'{a}n et al. (\citeyear{Her}) and \citet{plss} for ideas about dealing with censoring.

We denote the probability space by $(\Omega,\mathcal{F},P)$.
The covariate process describes the course of the patient
characteristics, for example, the course of the blood pressure and the white
blood cell count. We assume that a realization of this covariate
process is a function from $[0,\tau]$ to $\mathbb{R}^d$,
and that such a \textit{sample path} is continuous from the right with
limits from the left (cadlag). The covariates which \textit{must} be
included are those which both (i)~influence a doctor's treatment
decisions \textit{and} (ii)~possibly predict a patient's prognosis with
respect to the outcome of interest.  If such covariates would not be
observed, the assumption of no unmeasured confounding, mentioned in the
introduction, will not hold.

For the moment consider one single patient. We write $Z(t)$ for the
covariate- \textit{and} treatment values at time $t$. We assume that
$Z(t)$ takes values in $\mathbb{R}^m$, and that
$Z(t):\Omega\rightarrow \mathbb{R}^m$ is measurable for each
$t\in[0,\tau]$. Moreover, we assume that $Z$, seen as a
function on $[0,\tau]$, is cadlag. We write
$\overline{Z}_t=(Z(s)\dvtx 0\leq s\leq t)$ for the
covariate- and treatment history until time $t$, and $\overline{\mathcal{Z}}_t$
for the space of cadlag functions from $[0,t]$
to $\mathbb{R}^m$ in which $\overline{Z}_t$ takes it values.
Similarly, we write $\overline{Z}$ for the whole covariate- and
treatment history of the patient on the interval
$[0,\tau]$, and $\overline{\mathcal{Z}}$ for the space in
which $\overline{Z}$ takes its values. In this article we choose the
projection $\sigma$-algebra as the $\sigma$-algebra on
$\overline{\mathcal{Z}}_t$ and $\overline{\mathcal{Z}}$; measurability of
$Z(s)$ for each $s\leq t$ is then equivalent with measurability of the
random variable $\overline{Z}_t$.  For technical reasons, we include in
$Z$ a counter of the number of jump times of the measured treatment-
and covariate process. We suppose that observations on different
patients are independent.

\section{Counterfactual outcomes}

\begin{figure}[t]

\includegraphics{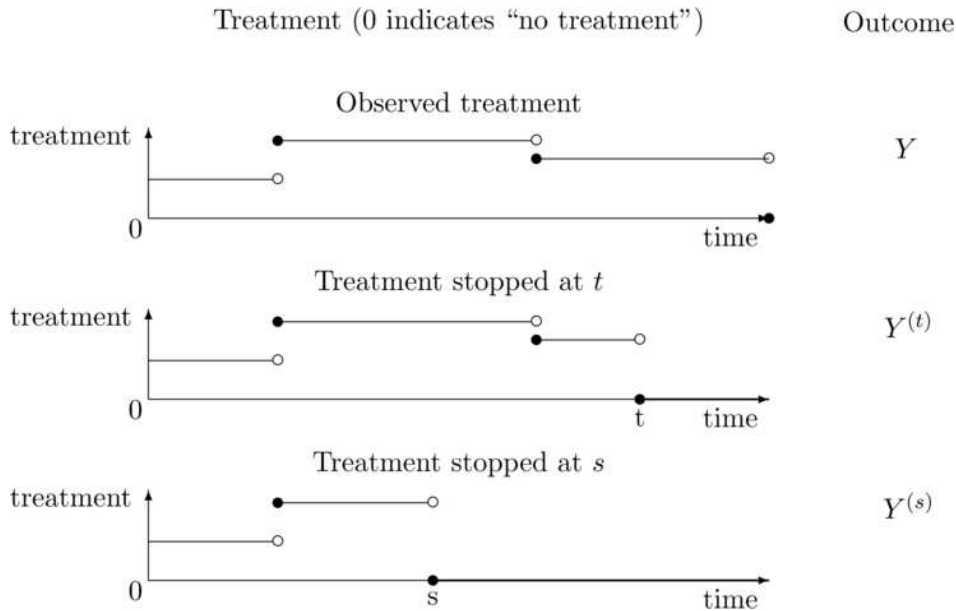}

\caption{Observed and counterfactual outcomes.}
\label{figuur1}
\end{figure}

Structural nested models are models for relations between so-called
counterfactuals. Consider for a moment just one patient. In reality
this patient received a certain treatment and had final outcome $Y$. If
his or her actual treatment had been stopped at time t, the patient's
final outcome would possibly have been different. The outcome he or
she would have had in that case we call $Y^{(t)}$. Of course,
$Y^{(t)}$ is generally not observed, because the patient's actual
treatment after $t$ is usually different from no treatment; it is a
counterfactual outcome. Instead of stopping treatment, one can also
consider switching to some kind of baseline treatment, for example, standard
treatment. Figure~\ref{figuur1} illustrates the nature of
counterfactual outcomes. We suppose that all counterfactual outcomes
$Y^{(t)}$, for $t\in[0,\tau]$ and for all patients, are
random variables on the probability space $(\Omega,\mathcal{F},P)$.

\section{No unmeasured confounding}
\label{nounsec}


To formalize the assumption of no unmeasured confounding, consider the
history of a particular patient. Decisions of the doctors at time $t$
may be based, in part, on recorded information on the state of the
patient and treatment before $t$, that is, on
$\overline{Z}_{t-}=(Z(s):0\leq s<t)$, but not on other
features predicting the outcome of the patient. In particular, given
$\overline{Z}_{t-}$, changes of treatment at time $t$ should be independent
of $Y^{(t)}$, the outcome of the patient in case he or she would not
have been treated after time $t$, given $\overline{Z}_{t-}$.

Note that $Y^{(t)}$ is an indication of the prognosis of the patient
which does not depend on treatment decisions at or after time $t$,
since it is the counterfactual outcome which we would have observed
if treatment would have been stopped at time~$t$. Only if treatment would have
no effect, the observed outcome $Y$ could play this role. This is why
Robins' assumption of no unmeasured confounding demands the
independence, given $\overline{Z}_{t-}$, of treatment decisions at time $t$
and $Y^{(t)}$. Similar conditions, though without time-dependence, can be
found in, for example, \citet{prop}.

The statement ``changes of treatment at $t$ should be independent of
$Y^{(t)}$, the outcome of the patient in case he or she would not have
been treated after time $t$, given $\overline{Z}_{t-}$'' is not a formal
statement: it includes conditioning of null events (since the
probability that treatment changes at $t$ may be $0$ for every fixed
$t$) on null events ($\overline{Z}_{t-}$).

To overcome this difficulty, we assume that the treatment process can
be represented by or generates a (possibly multivariate) counting
process $N$. For instance, $N(t)$ registers the number of
treatment changes until time $t$ and/or the number of times treatment
reached a certain level until time $t$. A counting process constructed
this way may serve as $N$ in the following. More about counting
processes can be found in, for example, Andersen et al. (\citeyear{ABGK}). We
assume that the treatment process $N$ has an intensity process.
Formally, such an intensity process $\lambda(t)$ is a
predictable process
such that $N(t)-\int_0^t\lambda(s)\,ds$ is a martingale.
The intensity $\lambda(t)$ with respect to
$\sigma(\overline{Z}_t)$ can be interpreted as the rate at which the
counting process $N$ jumps given the past treatment- and covariate
history $\overline{Z}_{t-}$.
\begin{assu}[(Bounded intensity process)]\label{int}
$N$ has an intensity process $\lambda(t)$ on
$[0,\tau]$ with respect to the filtration
$\sigma(\overline{Z}_t)$. This intensity process satisfies
the following regularity conditions:
\begin{longlist}[(a)]
\item[(a)] $\lambda$ is bounded by a constant which does not depend on
$\omega$,
\item[(b)] $\lambda(t)$ is continuous from the left.
\end{longlist}
\end{assu}

According to this assumption,
\begin{equation}
M(t)=N(t)-\int_0^t\lambda(s)\,ds\label{mg}
\end{equation}
is a
martingale on $[0,\tau]$ with respect to the filtration
$\sigma(\overline{Z}_t)$. Since most counting process
martingale theory deals with filtrations $\mathcal{F}_t$ which satisfy
the usual conditions ($\mathcal{F}_0$ contains all null sets and $\mathcal{F}_t=\bigcap_{s>t}
\mathcal{F}_s$), we mention that, under
\mbox{Assumption~\ref{int}}, $M(t)$ is also a
martingale with respect to $\sigma(\overline{Z}_t)^a$, the
usual augmentation of $\sigma(\overline{Z}_t)$. This
follows from Lemma~67.10 in \citet{RW}, since $M$ is cadlag.

Often, $N$ will be chosen to count the number of events of a certain
type concerning the treatment process (e.g., the number of
times treatment changed). At~$\tau$, the time the study ends,
treatment is stopped or switched to baseline treatment, so a natural
choice of $N$ will often jump at $\tau$ with positive probability.
However, jumps of $N$ at $\tau$ are not useful for estimation, and we
wish to avoid modeling jumps of $N$ at $\tau$. Therefore, we assume
that, with probability $1$, $N$ does not jump at $\tau$, and if a
natural choice of $N$ \textit{does} jump at $\tau$ with positive
probability, then we just adapt it, only at $\tau$, so that it does not
jump there.

We also make the following assumption.
\begin{assu}[(\textit{$Y^{\bolds{(\cdot)}}$ cadlag})]\label{Ycadlag} $Y^{\bolds{(\cdot)}}$ is a cadlag process.
\end{assu}

Within this framework, the assumption of no unmeasured confounding
could be operationalized as follows. The rate at which the counting
process $N$ jumps given past treatment- and covariate history is also
the rate at which $N$ jumps given past treatment- and covariate
history and $\overline{Y}^{(t)}=(Y^{(s)}\dvtx s\leq t)$. That
is, the following:
\begin{assu}[(No unmeasured confounding---formalization)]\label{intconf}
 The intensity process $\lambda(t)$ of $N$ with
respect to $\sigma(\overline{Z}_t)$ is also an intensity process of
$N$ with respect to $\sigma(\overline{Z}_t,\overline{Y}^{(t)})$.
\end{assu}

This can be interpreted as conditional independence (given
$\overline{Z}_{t-}$) of treatment decisions at time $t$ and
$(Y^{(s)}\dvtx s\leq t)$. This assumption is stronger than just
conditional independence of treatment decisions at $t$ and $Y^{(t)}$
as assumed in Robins (\citeyear{R92,Enc}). However, also $Y^{(s)}$ for $s<t$ is an
indication of the patient's prognosis upon which treatment decisions
at time $t$ ($>s$) should not depend.  Assumption~\ref{intconf} allows
us to use the usual counting processes framework. Under
Assumption~\ref{intconf},
$M(t)=N(t)-\int_{[0,t]}\lambda(s)\, ds$
is a martingale also with respect to $\sigma(\overline{Z}_t,\overline{Y}^{(t)})$ and
its usual augmentation $\sigma(\overline{Z}_t,\overline{Y}^{(t)})^a$.

This formalization of the assumption of no unmeasured confounding in
terms of compensators with respect to the filtration
$\sigma(\overline{Z}_t,\overline{Y}^{(t)})$ is a novel
feature of this paper relative to the previous literature on
structural nested models. Robins et al. (\citeyear{Aids}), \citet{Enc} and \citet{Kei} use
a Cox model for initiation and/or changes in treatment. However, none
of them formalized the assumption of no unmeasured confounders in
terms of compensators with respect to
$\sigma(\overline{Z}_t,\overline{Y}^{(t)})$. As a
consequence, they could not use the extensive theory on counting
process martingales to show the asymptotics of their estimators, which
then remained without proof.

\section{The model for treatment effect}
\label{model}

Structural nested models in continuous time model distributional
relations between $Y^{(t)}$ and $Y^{(t+h)}$, for $h>0$ small, through a
so-called infinitesimal shift-function $D$. Write $F$ for the
cumulative distribution function and $F^{-1}\dvtx (0,1)\mapsto
{\mathbb R}$ for its generalized inverse
\[
F^{-1}(p)=\inf\{x\dvtx F(x)\geq p\}.
\]
Then the infinitesimal shift-function $D$ is defined as
\begin{equation} \label{defd}
D(y,t;\overline{Z}_t)=\frac{\partial}{\partial h}\bigg|_{h=0}
\bigl(F_{Y^{(t+h)}| \overline{Z}_t}^{-1}\circ
F_{Y^{(t)}|\overline{Z}_t}\bigr)(y),
\end{equation}
the (right-hand) derivative of the quantile--quantile
transform which moves quantiles of the distribution of $Y^{(t)}$ to
quantiles of the distribution of $Y^{(t+h)}$ ($h\geq 0$), given
the covariate- and treatment history until time
$t$, $\overline{Z}_t$. Notice that for differentiability of $F_{Y^{(t+h)}|
\overline{Z}_t}^{-1}$ with respect to $h$ we need that a short duration of
treatment has only a small effect on the distribution of the outcome of
interest (see page~3, second paragraph), since $\lim_{h\downarrow 0}F_{Y^{(t+h)}|
\overline{Z}_t}$ must  be equal to $F_{Y^{(t)}| \overline{Z}_t}$.
\begin{exa}[(\textit{Survival of AIDS patients})]  \label{DAexa}
Robins, Blevins, Ritter and Wulfsohn (\citeyear{Aids}) describe an AIDS clinical trial to study the effect of
AZT treatment on survival in HIV-infected subjects. Embedded within
this trial was an essentially uncontrolled observational study of the
effect of prophylaxis therapy for PCP on survival. PCP, pneumocystis
carinii pneumonia, is an opportunistic infection that afflicts AIDS
patients. The aim of Robins et al. (\citeyear{Aids}) was to study the effect of this
prophylaxis therapy on survival. Thus, the outcome of interest $Y$ is
the survival time, and the treatment under study is prophylaxis for
PCP.  Once treatment with prophylaxis for PCP started, it was never
stopped.

Suppose that
\begin{equation}
D(y,t;\overline{Z}_t)=
(1-e^{\psi})1_{\{\mathrm{treated}\ \mathrm{at}\ t\}}.\label{Dvb1}\end{equation}
Then (see
Section~\ref{cmim} for details), for $t<Y$, withholding
treatment from $t$ onward leads to (with $\sim$ meaning ``is
distributed as'')
\begin{eqnarray} \label{DAsim}
Y^{(t)}-t &\sim& \int_{t}^Y e^{\psi
1_{\{\mathrm{treated}\ \mathrm{at}\ s\}}}\, ds\nonumber\\[-8pt]\\[-8pt]
&=&e^{\psi}\cdot \mathit{DUR}(t,Y)+1\cdot \bigl(Y-t-\mathit{DUR}(t,Y)\bigr)\qquad
\mbox{given }\overline{Z}_t,\nonumber
\end{eqnarray}
with $\mathit{DUR}(t,u)$ the duration of treatment in the interval $(t,u)$.
Thus, treated residual survival time ($t$ until $Y$) is multiplied by
$e^{\psi}$ by withholding treatment; compare this with accelerated failure
time models, see, for example, \citet{Cox}. This multiplication factor
$e^{\psi}$ should be interpreted in a distributional way. One of the
models studied in Robins et al. (\citeyear{Aids}) assumes that (\ref{DAsim}) is true even
with $\sim$ replaced by $=$ (though only for $t=0$). Notice that
supposing (\ref{DAsim}) to be true with $\sim$ replaced by $=$
would be much stronger.
\end{exa}
\begin{exa}[(\textit{Survival of AIDS patients})]
\label{DAexa2} Consider the situation from Example~\ref{DAexa} again.
In another model mentioned in Robins, Blevins, Ritter and Wulfsohn (\citeyear{Aids}) the factor with which
treated residual survival time is multiplied when treatment is
withheld can depend on the AZT treatment the patient received and
whether or not the patient had a history of PCP prior to start of PCP
prophylaxis. Since this was a clinical trial for AZT treatment, the
AZT treatment was described by a single variable $I_{\mathrm{AZT}}$ indicating
the treatment arm the patient was randomized to ($I_{\mathrm{AZT}}$ is $0$ or
$1$).  Whether or not the patient had a PCP history prior to start of
prophylaxis is described by an indicator variable $P(t)$.
$P(t)$ equals $1$ if the patient had PCP before or at $t$
\textit{and} before prophylaxis treatment started; otherwise
$P(t)$ equals $0$.  If
\[
D_{\psi_1,\psi_2,\psi_3}(y,t;\overline{Z}_t)
=\bigl(1-e^{\psi_1+\psi_2 P(t) +\psi_3 I_{\mathrm{AZT}}}\bigr)1_{\{\mathrm{treated}\ \mathrm{at}\
t\}},
\]
then (see
Section~\ref{cmim} for details) withholding prophylaxis treatment
from $t$ onward leads to
\begin{equation} \label{DAsim2}
Y^{(t)}-t \sim \int_{t}^Y e^{1_{\{\mathrm{treated}\,\mathrm{at}\,s\}}(\psi_1+\psi_2 P(s) +\psi_3 I_{\mathrm{AZT}})}\, ds\qquad
\mbox{given }\overline{Z}_t,\end{equation}
for $t<Y$.
\end{exa}
\begin{exa}[(\textit{Incorporating a-priori biological knowledge})] \label{prior}
Following \citet{Enc},
again consider survival as the outcome of
interest.
Suppose that it is known that treatment received at time $t$
only affects survival
for patients destined to die by time $t+5$ if they would receive
no further treatment. An example would be a setting
in which failure is death from an infectious disease, the treatment
is a preventive antibiotic treatment which is of no benefit
unless the subject is already infected and, if death occurs, it
always does within five weeks from the time of initial unrecorded
subclinical infection.
In that case, as
remarked in \citet{Enc}, the natural restriction on $D$ is that
\[
D(y,t;\overline{Z}_t)=0\qquad \mathrm{if}\ y-t>5.
\]
\end{exa}

More biostatistical examples of models for $D$ can be found in,
for example, \citet{smoke}, Witteman et al. (\citeyear{hyp}), \citet{Enc} and Keiding et al.
(\citeyear{Keiding}).

\begin{figure}[t]

\includegraphics{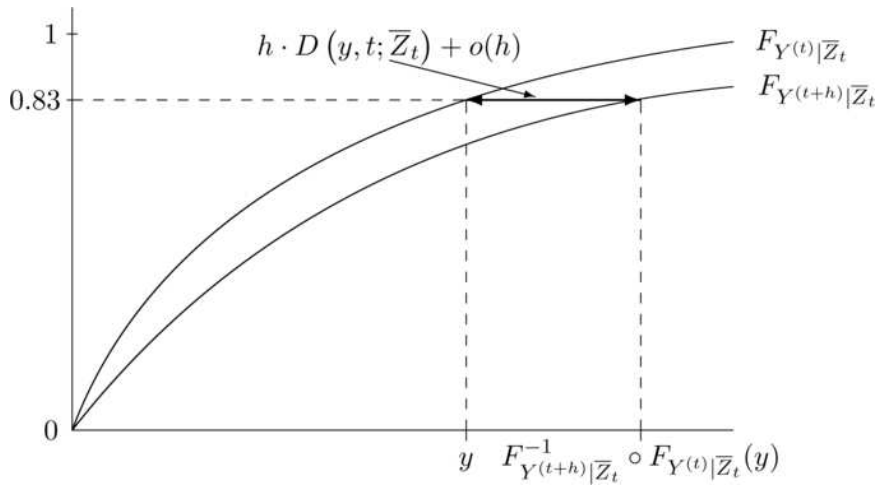}

\caption{Illustration of the infinitesimal shift-function $D$.}
\label{Dfig}
\end{figure}

$D(y,t;\overline{Z}_t)$ can be interpreted as the infinitesimal effect
on the outcome $Y$ of the treatment actually given in
the time-interval $[t,t+h)$ (relative to baseline
treatment). To be
more precise, from the definition of $D$ we have
\[
h\cdot D(y,t;\overline{Z}_t)=
\bigl(F_{Y^{(t+h)}| \overline{Z}_t}^{-1}\circ
F_{Y^{(t)}|\overline{Z}_t}\bigr)(y)-y+o(h).
\]
In Figure~\ref{Dfig} this is sketched. $y$ in the picture is the
$0.83$th quantile of the distribution of $Y^{(t)}$ given $\overline{Z}_t$. For
$h>0$, the $0.83$th quantile of the distribution of $Y^{(t+h)}$ given $\overline{Z}_t$
is $y+h\cdot D(y,t;\overline{Z}_t) +o(h)$. Thus, to shift
from quantiles of the distribution of $Y^{(t)}$ to the distribution of
$Y^{(t+h)}$ given $\overline{Z}_t$ ($h>0$) is approximately the same as to just add
$h\cdot D(y,t;\overline{Z}_t)$ to those quantiles.  For example, if
$F_{Y^{(t+h)}|\overline{Z}_t}$ lies to the right of $F_{Y^{(t)}|\overline{Z}_t}$ for $h>0$, then
treatment between $t$ and $t+h$ increases the outcome (in
distribution), and $D(\cdot,t;\overline{Z}_t)$ is greater than $0$.

Consider again this interpretation of $D$ as the infinitesimal effect
of treatment given in $[t,t+h)$. If the outcome of interest
is survival, then $D(y,t;\overline{Z}_t)$ should be zero if
$\overline{Z}_t$ indicates the patient is dead at time $t$.  Indeed, in
that case $F_{Y^{(t+h)}|\overline{Z}_t}$ and
$F_{Y^{(t)}|\overline{Z}_t}$ should be almost surely the same for
every $h\geq 0$, since withholding treatment after death does not
change the survival time. Thus, $F_{Y^{(t+h)}|\overline{Z}_t}^{-1}\circ
F_{Y^{(t)}|\overline{Z}_t}(y)$ is constant in $h$ for
$h\geq 0$ and, therefore, $D(y,t;\overline{Z}_t)=0$. However,
this reasoning is not precise because of the complication of null
sets. We will therefore just formally define
$D(y,t;\overline{Z}_t)$ to be zero if the outcome of
interest is survival and $\overline{Z}_t$ indicates the patient is
dead at time $t$.

It can be shown that $D\equiv 0$ if treatment does not affect the
outcome of interest, as was conjectured in \citet{Enc}. To be more
precise, \citet{Lok} shows that, for example, $D\equiv 0$ if and only
if, for every $h>0$ and $t$, $Y^{(t+h)}$ has the same distribution as
$Y^{(t)}$ given $\overline{Z}_t$. That is, $D\equiv 0$ if and only
if ``at any time $t$, whatever patient characteristics are selected at
that time ($\overline{Z}_t$), stopping `treatment as given' at some fixed time
after $t$ would not change the distribution of the outcome in patients
with these patient characteristics.''

In the rest of this article $D_\psi$ will always indicate a correctly
specified parametric model for $D$, with $D\equiv 0$ if $\psi=0$.

\section{Mimicking counterfactual outcomes}
\label{cmim}

\begin{figure}[b]

\includegraphics{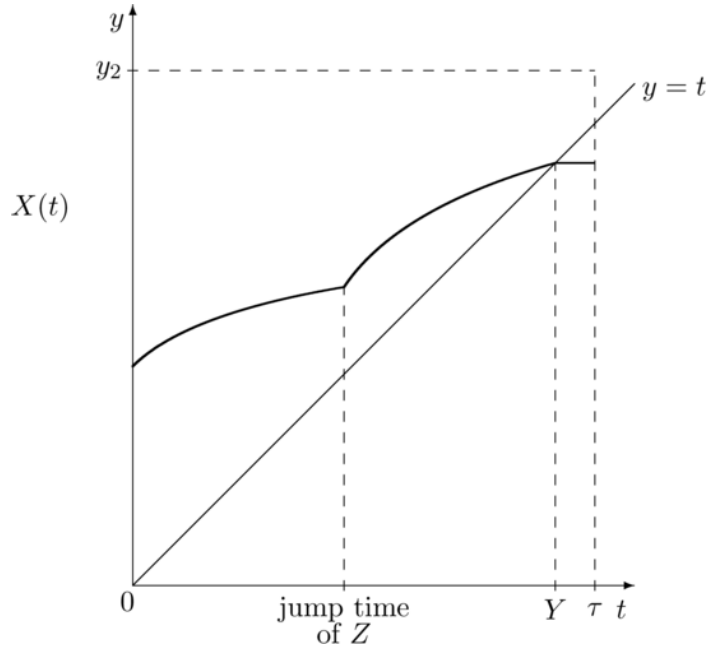}

\caption{An example of a solution $X(t)$ to the differential
equation $X'(t)=D(X(t),t;\overline{Z}_t)$ with final condition
$X(\tau)=Y$ in case the outcome is survival time.}
\label{Dtfig}
\end{figure}

Define $X(t)$ as the
continuous solution to the differential equation
\begin{equation} \label{Xdef}
X'(t)=D(X(t),t;\overline{Z}_t)
\end{equation}
with final condition $X(\tau)=Y$, the observed outcome (see
Figure~\ref{Dtfig}). Then $X(t)$ mimics $Y^{(t)}$ in the
sense that it has the same distribution as $Y^{(t)}$ given
$\overline{Z}_t$. This rather surprising result was conjectured in
\citet{Enc} and proved in Lok (\citeyear{Lok,Mimarx}). To prove this result, we need
the following consistency assumption.
\begin{assu}[(Consistency)]\label{inst}
$Y^{(\tau)}$ has the same distribution as $Y$ given~$\overline{Z}_{\tau}$.
\end{assu}

Notice that since by assumption no treatment was given after time
$\tau$ and since treatment is right-continuous, there is no difference
in treatment between $Y^{(\tau)}$ and~$Y$. We suppose that this
assumption holds, and we also suppose that a short duration of
treatment has only a small effect on the distribution of the outcome
of interest ($\lim_{h\downarrow
  0}F_{Y^{(t+h)}|\overline{Z}_t}(y)\rightarrow
F_{Y^{(t)}|\overline{Z}_t}(y)$). Under these assumptions and
regularity conditions only, Lok (\citeyear{Lok,Mimarx}) proved that indeed
equation (\ref{Xdef}) has a unique solution for every
$\omega\in\Omega$, and that this solution $X(t)$ mimics $Y^{(t)}$ in
the sense that $X(t)$ has the same distribution as $Y^{(t)}$ given
$\overline{Z}_t$ (see Appendix~\ref{XTap}). Throughout this article we
will assume that this result holds true.
\begin{exa} \label{DAexai}
Survival of AIDS patients (continuation of Example~\textup{\ref{DAexa}}).
Suppose that
\[D(y,t;\overline{Z}_t)=
(1-e^{\psi}) 1_{\{\mathrm{treated}\ \mathrm{at}\ t\}}.\]
Then
\[
X(t)=t+ \int_{t}^Y e^{\psi
1_{\{\mathrm{treated}\ \mathrm{at}\ s\}}}\, ds
\]
if $Y>t$, and $X(t)=Y$ for $t\geq Y$.
\end{exa}

Suppose now that one has a correctly specified parametric model for
the infinitesimal shift-function $D$, $D_\psi$. Then one can calculate
``$X_\psi(t)$,'' the solution to
\begin{equation} \label{Xpsi}
X_\psi'(t)=
D_{\psi}(X_\psi(t),t;\overline{Z}_t)
\end{equation}
with final condition $X_\psi(\tau)=Y$.  For the true $\psi$,
$X_\psi(t)$ has the same distribution as $Y^{(t)}$, the outcome with
treatment stopped at $t$, even given all patient-information at time
$t$, $\overline{Z}_t$. So instead of the unobservable $Y^{(t)}$'s, we
have the observable $X_{\psi}(t)$'s which for the true $\psi$ mimic
the $Y^{(t)}$'s. Although we do not know the true $\psi$, this result
turns out to be very useful, both for estimating $\psi$
(Sections~\ref{ees} and~\ref{asnorm}) and for testing
[Section~\ref{test}; notice that when testing whether treatment
affects the outcome (i.e., whether $D\equiv 0$), $X$ can simply be calculated
from the data ($X\equiv Y$) under the null hypothesis of no treatment
effect].

\section{Local rank preservation}
\label{loc}

Previous applications of structural nested models [see, e.g., Robins et al. (\citeyear{Aids}), \citeauthor{smoke}
(\citeyear{smoke}), Witteman et al. (\citeyear{hyp}) and
Keiding et al. (\citeyear{Keiding})] have assumed the so-called
local rank preservation condition.
Local rank preservation states that $Y^{(t)}$ is a local solution to
(\ref{Xdef}). However, if $Y^{(t)}$ is locally a solution to
(\ref{Xdef}), it is usually also globally a solution to (\ref{Xdef});
see, for example, Theorem~\ref{isun} in the \hyperref[difeqap]{Appendix}. Hence, if one knew the
parameter $\psi$, every $Y^{(t)}$ would be a deterministic function of
the observed data. Deterministic dependence of counterfactuals on the
observed data is a very strong condition, which, though untestable, is
generally considered implausible. The previous literature
[see, e.g., \citeauthor{Enc} (\citeyear{Enc}), Robins et al. (\citeyear{Aids}),
\citeauthor{smoke} (\citeyear{smoke}) and Keiding et al.
(\citeyear{Keiding})] acknowledged this problem, and conjectured that
the assumption of local rank preservation could be relaxed in
continuous time [since it is known that the assumption of local rank
preservation can be relaxed for structural nested models in discrete
time; this was pointed out by \citeauthor{RobW} (\citeyear{RobW}), and
Lok et al. (\citeyear{SNart}) provided a proof]. See \citet{Enc}
for a more elaborate discussion. The following example describes a
setting where the assumption of local rank preservation is
implausible.
\begin{exa}[(\textit{Survival of AIDS patients and local rank
preservation})]
In the situation of Example~\ref{DAexa}, consider the following
thought experiment. Suppose that two patients had the same covariate
history until time $t$, and both received the same constant treatment
to prevent PCP until time $t$
(equal $\overline{Z}_t$). Suppose, furthermore, that both patients
received no treatment after time $t$, that they did not have PCP
before time $t$ and both died at the same time $u>t$ (for both,
$Y^{(t)}=u$). Possibly, the first patient would have had PCP at some
time $s<t$ and would have died from it before $u$ in case he or she
would not have been treated. Possibly, the other patient would not
have had PCP in case he or she would not have been treated, and would
have died at the same time $u>t$ as without treatment. Thus, it is
easy to imagine that these patients would have had different outcomes
under no treatment (different $Y^{(0)}$). However, the
assumption of local rank preservation excludes this possibility.
\end{exa}

Local rank preservation is a very strong condition, for which
structural nested models have previously been attacked. In fact, this
article shows that the assumption of local rank preservation is not
needed for structural nested models. However, proofs would be much
easier under rank preservation; for details, see the remarks before the
proofs of Theorems~\ref{thec} and~\ref{mgc}. See also \citet{Enc} for a
more informal reasoning.

\section{Estimation of treatment effect}
\label{ees}

To estimate the infinitesimal shift-\break function $D$, \citet{Enc} proposes
to use a (semi-)parametric model to predict future treatment ($N$ in our
case) on the basis of past treatment- and covariate history
$\overline{Z}_{t-}$.  This may seem odd, since prediction of treatment
is not what we are interested in. However, we will show that such
a model to predict treatment changes can indeed be a tool to get
unbiased estimating equations for the parameter $\psi$ in the model
for $D$.  Moreover, often doctors may have a better understanding, at
least qualitatively, about how decisions about treatment were made
than about the effect of the treatment. In what follows we will assume
that $\lambda_\theta$ is a correctly specified parametric model for
the intensity $\lambda$ of $N$.

Recall from Section~\ref{nounsec} that, under no unmeasured confounding
(Assumption~\ref{intconf}), $Y^{(t)}$ does not contain information
about treatment changes given past treatment- and covariate history
$\overline{Z}_{t-}$. Since $X(t)$ has the same distribution as
$Y^{(t)}$ given $\overline{Z}_t$, one could expect that also $X(t)$
does not contain information about treatment changes given
$\overline{Z}_{t-}$. Unfortunately, this reasoning is not precise: we
have to somehow deal with null sets since the probability that
treatment changes at $t$ given past covariate- and treatment history
is often equal to $0$ for each $t$. In Section~\ref{mgs} we will show
how this can be dealt with.

In the current section we present a class of unbiased estimating
equations for $\theta$ and $\psi$. These will be used for the proof in
the next section, but they are also of interest in their own. In
Section~\ref{mgs} we will see that these estimating equations are in
fact martingales, for the true parameters $\theta_0$ and $\psi_0$.

Recall from Section~\ref{nounsec} that, under no unmeasured confounding, we have the martingale
$M(t)=N(t)-\int_{[0,t]}\lambda(s)\,ds$ with respect to the
filtration $\sigma(\overline{Z}_t,\overline{Y}^{(t)})$ and its usual
augmentation. From this martingale we can construct a whole family of
martingales. If $h_t(Y^{(t-)},\overline{Z}_{t-})$
is a $\sigma(\overline{Z}_t,\overline{Y}^{(t)})^a$-predictable process, then, under
regularity conditions,
\[
\int_0^t h\,dM
=\int_0^t h_s\bigl(Y^{(s-)},\overline{Z}_{s-}\bigr)\, dM(s)
\]
is a martingale with respect to $\sigma(\overline{Z}_t,\overline{Y}^{(t)})^a$.
For a more formal statement, we first make sure that
$h_t(Y^{(t-)},\overline{Z}_{t-})$ is predictable. We
put the following restriction on the functions $h_t$ we consider here:
\begin{rest}\label{hreg}
When in this section we consider functions $h_t$ from
$\mathbb{R}\times \overline{\mathcal{Z}}_{t-}$, we assume that they are
measurable and satisfy the following:
\begin{longlist}[(a)]
\item[(a)] $h_t$ is bounded by a constant which does not depend on $t$ and
$\overline{Z}$,
\item[(b)] for all $t_0\in[0,\tau]$, $y_0\in\mathbb{R}$ and
$\omega\in\Omega$,
$h_t(y,\overline{Z}_{t-}(\omega))\rightarrow
h_{t_0}(y_0,\overline{Z}_{t_0-}(\omega))$ when
$y\rightarrow y_0$ and $t\uparrow t_0$.
\end{longlist}
\end{rest}

For such $h_t$, $h_t(Y^{(t-)},\overline{Z}_{t-})$ is a
$\sigma(\overline{Z}_t,\overline{Y}^{(t)})^a$-predictable process:
\begin{lem} \label{hpred}
Suppose that $Y^{\bolds{(\cdot)}}$ is cadlag \textup{(}Assumption~\textup{\ref{Ycadlag})}. Then
$h_t(Y^{(t-)},\break\overline{Z}_{t-})$ is a
$\sigma(\overline{Z}_t,\overline{Y}^{(t)})^a$-predictable process for any $h_t$
satisfying Restriction~\textup{\ref{hreg}}.
\end{lem}
\begin{pf*}{Proof} $h_t(Y^{(t-)},\overline{Z}_{t-})$ is adapted. It is also
left-continuous: because $Y^{\bolds{(\cdot)}}$ is cadlag,
$\lim_{t\uparrow t_0} Y^{(t-)}=Y^{(t_0-)}$ exists.
\end{pf*}

 Thus, we come to the following lemma:
\begin{lem} \label{lemec} Under Assumptions~\textup{\ref{int}
(}bounded intensity process\textup{), \ref{Ycadlag}} \textup{(}$Y^{\bolds{(\cdot)}}$ cadlag\textup{)}
and~\ref{intconf} (no unmeasured confounding),
\[
\int_0^t h_s\bigl(Y^{(s-)},\overline{Z}_{s-}\bigr)
\bigl(dN(s)-\lambda(s)\,ds\bigr)
\]
is a martingale on $[0,\tau]$ with respect to
$\sigma(\overline{Z}_t,\overline{Y}^{(t)})^a$ for all $h_t$ satisfying
Restriction~\textup{\ref{hreg}}.
\end{lem}
\begin{pf*}{Proof}  $M(t)=N(t)-\int_0^t\lambda(s)\,ds$ is
a martingale on $[0,\tau]$ with respect to
$\sigma(\overline{Z}_t,\overline{Y}^{(t)})^a$, because of Assumption~\ref{intconf}. It
is of integrable variation\break [$E\int_0^t | dM(s)| \leq
E\int_0^t dN(s)+\lambda(s)\, ds = 2 E \int_0^t \lambda(s)\, ds$, and
$\lambda$ is bounded (Assumption~\ref{int})]. Because of
Lemma~\ref{hpred},
$h(t)=h_t(Y^{(t-)},\overline{Z}_{t-})$ is a
$\sigma(\overline{Z}_t,\overline{Y}^{(t)})^a$-predictable process. It is also
bounded [Restriction~\ref{hreg}(a)]. Hence,\break
$\int_0^t h(s)\,dM(s)=
\int_0^t h_s(Y^{(s-)},\overline{Z}_{s-})
(dN(s)-\lambda(s)\,ds)$
is an integral of a bounded predictable process with respect to a
martingale of integrable variation, and, therefore, a~$\sigma(\overline{Z}_t,\overline{Y}^{(t)})^a$-martingale.
\end{pf*}

To construct unbiased estimating equations for
$(\theta_0,\psi_0)$, we need to assume that the probability
that $N(\cdot)$ and $Y^{\bolds{(\cdot)}}$ jump at the same time is zero. This assumption
is a formalization of the assumption of no instantaneous treatment
effect as proposed in \citet{Enc}, which can be seen as follows.  Given
$\overline{Z}_{t-}$ and $\overline{Y}^{(t-)}$, $N$ jumps at
$t$ with rate $\lambda(t)$ (Assumption~\ref{intconf}, no unmeasured
confounding).  $Y^{\bolds{(\cdot)}}$ is a cadlag process
(Assumption~\ref{Ycadlag}), which thus for every $\omega\in\Omega$
jumps at most countably many times on the finite time interval
$[0,\tau]$. Therefore, if $Y^{\bolds{(\cdot)}}$ and $N$ would jump at the
same time with positive probability, this would imply a dependence of
these jumps; the obvious interpretation of this dependence would be
that a change of treatment instantaneously affects the outcome of
interest.
\begin{assu}[(\textit{No instantaneous
treatment effect})] \label{jump}  The probability that there exists a $t$ such that
$N(\cdot)$ and $Y^{\bolds{(\cdot)}}$ both jump at time $t$ is $0$.
\end{assu}

Notice that this excludes estimation of the effect of point exposures.
For example, if treatment is surgery or another point exposure given
at some time $t$, the outcome under ``treatment stopped at time $t$''
will typically jump at time $t$ if treatment affects the outcome of
interest, at the same time as the treatment itself. However, this
assumption does not exclude the possibility that the outcome differs
depending on whether a patient is treated or not at a certain point in
time.  For example, $Y^{(t+)}$ and $Y^{(t-)}$ may be different when a virus
is contacted at time $t$. The model in this article can accommodate
differences between $Y^{(t+)}$ and $Y^{(t-)}$, as long as the
probability that the observed treatment changes is $0$ at that precise
time. Or, in more generality, as long as the probability that $N$
jumps at the same time is $0$. This was previously noticed in
\citet{corr}, Section~8. The estimating procedures in this article do
not deal with instantaneous treatment effects.

Suppose that the above conditions hold and that
$(X(t),\overline{Z}_t)\sim (Y^{(t)},\overline{Z}_t)$ for
$t\in[0,\tau]$ (see Section~\ref{cmim}). Then if $D_{\psi}$
and $\lambda_\theta$ are correctly specified (parametric) models for
$D$ and $\lambda$, respectively, each choice of $h_t$ satisfying
Restriction~\ref{hreg} leads to an unbiased estimating equation for
both the parameter of interest $\psi$ and the (nuisance) parameter
$\theta$:
\begin{theorem} \label{thec}
Suppose that Assumptions~\textup{\ref{int} (}bounded intensity process\textup{),
\ref{Ycadlag} (}$Y^{\bolds{(\cdot)}}$ cadlag\textup{), \ref{intconf}} (no unmeasured
confounding) and~\textup{\ref{jump}} (no instantaneous treatment effect) are
satisfied. Suppose also that, for every $t\in[0,\tau]$,
$X(t)$ has the same distribution as $Y^{(t)}$ given $\overline{Z}_t$.  Then
\[
E\int_0^\tau  h_t(X(t),\overline{Z}_{t-})
 \bigl(dN(t)-\lambda(t)\,dt\bigr)=0
 \]
for each $h_t$ satisfying Restriction~\textup{\ref{hreg}}. Thus, if
$D_\psi$ and $\lambda_\theta$ are correctly specified parametric
models for $D$ and $\lambda$, respectively,
\[
P_n\int_0^\tau  h_t(X_{\psi}(t),\overline{Z}_{t-})
 \bigl(dN(t)-\lambda_{\theta}(t)\,dt\bigr)=0,
 \]
with $P_n$ the empirical measure $P_n \,X=1/n \sum_{i=1}^n X_i$, is an
unbiased estimating equation for $(\theta_0,\psi_{0})$, for
each $h_t$ satisfying Restriction~\textup{\ref{hreg}}. $h_t$ here is allowed to
depend on $\psi$ and $\theta$, as long as it satisfies
Restriction~\textup{\ref{hreg}} for $(\theta_0,\psi_{0})$.
\end{theorem}

As before, $X_{\psi}(t)$ here is the
continuous solution of (\ref{Xpsi}),
$X_\psi'(t)=D_\psi(X_\psi(t),t;\break\overline{Z}_t)$
with boundary condition $X_\psi(\tau)=Y$. Moreover, as
before, we assume that for all $D_\psi$ we have existence and
uniqueness of such solutions; Theorem~\ref{isun} in the appendix
provides sufficient conditions for that.

Under local rank preservation (see Section~\ref{loc}), $X(t)=Y^{(t)}$
for each $t$. In that case Theorem~\ref{thec} follows immediately
from Lemma~\ref{lemec}. However, as argued in Section~\ref{loc}, local
rank preservation is generally considered implausible.
\begin{pf*}{Proof of Theorem \ref{thec}} We have to show that
\[
\int_0^\tau  h_t(X(t),\overline{Z}_{t-})
\bigl(dN(t)-\lambda(t)\,dt\bigr)
\]
has expectation zero for all $h_t$ satisfying Restriction~\ref{hreg}.
To do that, we prove that it has the same expectation as
\[\int_0^\tau  h_t\bigl(Y^{(t-)},\overline{Z}_{t-}\bigr)
\bigl(dN(t)-\lambda(t)\,dt\bigr),\]
which has expectation zero because of Lemma~\ref{lemec}. We will first
show that the terms with $dN$ have the same expectation, that is,
\begin{equation}
E\Biggl(\sum_{t\leq \tau,\Delta
N(t)=1}h_t(X(t),\overline{Z}_{t-})\Biggr)=
E\Biggl(\sum_{t\leq \tau,\Delta
N(t)=1}h_t\bigl(Y^{(t-)},\overline{Z}_{t-}\bigr)\Biggr)\label{firste}.
\end{equation}
After that we show that the terms with $\lambda(t)\, dt$ have the same
expectation, that is,
\begin{equation} E\biggl(\int_0^\tau
h_t(X(t),\overline{Z}_{t-}) \lambda(t)\,dt
\biggr) =E\biggl(\int_0^\tau
h_t\bigl(Y^{(t-)},\overline{Z}_{t-}\bigr) \lambda(t)\,dt\biggr)
\label{seconde}.\end{equation}
As we will see below, (\ref{firste}) and (\ref{seconde}) have to be
proved separately, since we do not have or expect that
$(X(s),\overline{Z}_t)\sim
(Y^{(s)},\overline{Z}_t)$
for $s<t$; we only have this for $s\geq t$. Therefore, the
approximations below have to be chosen carefully.

At first we prove (\ref{firste}), by approximating these sums and
showing that the approximations have the same expectation. Next we
show that the approximations converge and that (\ref{firste}) follows
with Lebesgue's dominated convergence theorem.

Define $T_1=\inf\{t\dvtx N(t)=1\}$,
$T_2=\inf\{t\dvtx N(t)=2\}$, etc., the jump times of
the counting process $N$ in the interval $[0,\tau]$. They
are measurable [e.g., because of \citeauthor{RW} (\citeyear{RW}),
Lemma~74.4]. Note that the number of jumps in $[0,\tau]$ is
almost surely finite because $N$ is integrable (it has a bounded
intensity process). In the following read
$h_{T_j}(Y^{(T_j-)},\overline{Z}_{T_j-})=0$ if
there is no $j$th jump of $N$ in the interval $[0,\tau]$.

Next split up the interval $[0,\tau]$ in intervals of equal
length: for $K\in{\mathbb N}$ fixed, put $\tau_k=k\tau/K$,
$k=0,\ldots,K$.  Fix $K$ for the moment. The right-hand side of
equation~(\ref{firste}) is harder to approximate than the left-hand
side, both because $Y^{(t)}$ does not need to be continuous in $t$
while $X(t)$ does and because knowing $Y^{(t)}$ and $\overline{Z}_t$
does not imply knowing $Y^{(s)}$ for $s<t$ and we do not have or
expect $(X(s),\overline{Z}_t)\sim
(Y^{(s)},\overline{Z}_t)$ for $s<t$.  The
approximations we choose are
\begin{eqnarray}\label{apY}
\qquad \sum_{{\Delta N(t)=1, t\leq \tau}}
h_t\bigl(Y^{(t-)},\overline{Z}_{t-}\bigr)
&=&\sum_{j=1}^\infty
h_{T_j}\bigl(Y^{(T_{j}-)},\overline{Z}_{T_j-}\bigr)\nonumber\\[-8pt]\\[-8pt]
\qquad &\approx&
\sum_{j=1}^\infty\sum_{k=0}^{K-1}
1_{(\tau_k,\tau_{k+1}]}(T_j)
h_{\tau_{k}}\bigl(Y^{(\tau_{k+1})},\overline{Z}_{\tau_{k}-}\bigr)\nonumber
\end{eqnarray}
and
\begin{eqnarray}\label{apX}
\qquad \sum_{{\Delta N(t)=1,t\leq \tau}}h_t(X(t),\overline{Z}_{t-})
&=&\sum_{j=1}^\infty
h_{T_j}(X(T_{j}),\overline{Z}_{T_j-})\nonumber\\[-8pt]\\[-8pt]
\qquad &\approx& \sum_{j=1}^\infty\sum_{k=0}^{K-1}
1_{(\tau_k,\tau_{k+1}]}(T_j)
h_{\tau_{k}}\bigl(X(\tau_{k+1}),\overline{Z}_{\tau_{k}-}\bigr).\nonumber
\end{eqnarray}

To show that these approximations have the same expectation, we use
that
$(X(\tau_{k+1}),\overline{Z}_{\tau_{k+1}})\sim
(Y^{(\tau_{k+1})},\overline{Z}_{\tau_{k+1}})$. Therefore,
also
\[
1_{(\tau_k,\tau_{k+1}]}(T_j)
h_{\tau_{k}}(X(\tau_{k+1}),\overline{Z}_{\tau_{k-}})\sim
1_{(\tau_k,\tau_{k+1}]}(T_j)
h_{\tau_{k}}\bigl(Y^{(\tau_{k+1})},\overline{Z}_{\tau_{k-}}\bigr)
\]
[notice that $1_{(\tau_k,\tau_{k+1}]}(T_j)$ is a
function of $\overline{Z}_{\tau_{k+1}}$]. Hence, the expectation of each of
the terms on the right-hand side of (\ref{apY}) is equal to the expectation of the
corresponding term on the right-hand side of (\ref{apX}). Since $h_t$ is bounded
[Restriction~\ref{hreg}(a)] and the expected number of jump times $T_j$
is finite ($N$ is integrable), this implies that the expectation
of the right hand-side of equation~(\ref{apY}) is equal to the
expectation of the right-hand side of equation~(\ref{apX}).

Equation~(\ref{firste}) follows if the expectation of the
approximations in (\ref{apY}) and (\ref{apX}) converges to the
right-hand side and left-hand side of equation~(\ref{firste}),
respectively. This convergence is harder to show for (\ref{apY}) than
for (\ref{apX}), since $Y^{\bolds{(\cdot)}}$ may jump [while $X(\cdot)$ does not, by
construction].  Fix $j$ for a moment. Define $\tau_k^j$ and
$\tau_{k+1}^j$ as the grid points such that
$T_j\in(\tau_k^j,\tau_{k+1}^j]$. As $K\rightarrow\infty$,
$\tau^j_{k+1} \downarrow T_j$, so that since $Y^{(\bolds{\cdot})}$ is cadlag,
$Y^{(\tau_{k+1}^j)}\rightarrow Y^{(T_j)}$. Moreover, as $K\rightarrow
\infty$, $\tau_k^j\uparrow T_j$, so that because of
Restriction~\ref{hreg} on $h$,
$h_{\tau_k^j}(Y^{(\tau_{k+1}^j)},\overline{Z}_{{\tau_k^j}-})
\rightarrow h_{T_j}(Y^{(T_j)},\overline{Z}_{T_j-})$.
Combining this for all $j$ leads to
\[
\sum_{j=1}^\infty\sum_{k=0}^{K-1}
1_{(\tau_k,\tau_{k+1}]}(T_j)
h_{\tau_{k}}\bigl(Y^{(\tau_{k+1})},\overline{Z}_{\tau_{k-}}\bigr)
\rightarrow \sum_{j=1}^\infty
h_{T_j}\bigl(Y^{(T_{j})},\overline{Z}_{T_j-}\bigr)
\] as
$K\rightarrow \infty$ for every $\omega$ for which the number of jumps
of $N$ is finite, so for almost every $\omega\in\Omega$. $h_t$ is
bounded [Restriction~\ref{hreg}(a)] and the left-hand side is bounded by
the number of jumps of $N$ times this bound. The expectation of that
is finite because $N$ is integrable. Thus, Lebesgue's dominated
convergence theorem can be applied, and
\begin{equation} \label{naarplus}
\qquad E\Biggl(\sum_{j=1}^\infty\sum_{k=0}^{K-1}
1_{(\tau_k,\tau_{k+1}]}(T_j)
h_{\tau_{k}}\bigl(Y^{(\tau_{k+1})},\overline{Z}_{\tau_{k}-}\bigr)\Biggr)
\rightarrow E\Biggl(\sum_{j=1}^\infty
h_{T_{j}}\bigl(Y^{(T_{j})},\overline{Z}_{T_j-}\bigr)\Biggr)
\end{equation}
as
$K\rightarrow \infty$. Because with probability one $Y^{\bolds{(\cdot)}}$ and
$N$ do not jump at the same time
(Assumption~\ref{jump} of no instantaneous treatment effect),
\[
\sum_{j=1}^\infty h_{T_j}\bigl(Y^{(T_{j})},\overline{Z}_{T_j-}\bigr)
=\sum_{j=1}^\infty
h_{T_j}\bigl(Y^{(T_{j}-)},\overline{Z}_{T_j-}\bigr)\qquad\mathrm{a.s.},
\]
so that we can replace $Y^{(T_{j})}$ by
$Y^{(T_{j}-)}$ on the right-hand side of
(\ref{naarplus}). Therefore, indeed, the expectation of the approximation
in (\ref{apY}) converges to the expectation of the left-hand side of
(\ref{apY}).
The same reasoning shows this for (\ref{apX}).
Here less caution is necessary since $X(t)$ is continuous in $t$. That
concludes the proof of equation~(\ref{firste}).

Next we prove (\ref{seconde}), also by approximation. Here, too, we
show that the approximations have the same expectation and that
(\ref{seconde}) follows with Lebesgue's dominated convergence theorem.

Divide the interval $[0,\tau]$ as above.
The approximations we choose here are
\begin{equation} \label{apY2}
\int_0^\tau  h_{t}\bigl(Y^{(t-)},\overline{Z}_{t-}\bigr)
\lambda(t)\,dt
\approx \sum_{k=0}^{K-1}
h_{\tau_k}\bigl(Y^{(\tau_k)},\overline{Z}_{\tau_k-}\bigr)
\lambda(\tau_k)(\tau_{k+1}-\tau_k)\end{equation}
and
\begin{equation} \label{apX2}
\int_0^\tau h_{t}(X(t),\overline{Z}_{t-})
\lambda(t)\, dt \approx \sum_{k=0}^{K-1}
h_{\tau_k}(X(\tau_k),\overline{Z}_{\tau_k-})
\lambda(\tau_k)(\tau_{k+1}-\tau_k).
\end{equation}
Because $(X(\tau_k),\overline{Z}_{\tau_k})\sim
(Y^{(\tau_k)},\overline{Z}_{\tau_k})$ and
$\lambda(\tau_k)$ are a measurable function of
$\overline{Z}_{\tau_k}$ (Assumption~\ref{int}, bounded intensity
process), the expectation of each of the terms in (\ref{apY2}) is
equal to the expectation of the corresponding term in
(\ref{apX2}). Thus, the expectations of these approximations are equal.

Equation~(\ref{seconde}) follows if the expectation of the
approximations in~(\ref{apY2}) and~(\ref{apX2}) converge to the
right-hand side and left-hand side of equation~(\ref{seconde}),
respectively. This convergence is also harder to show for (\ref{apY2})
than for (\ref{apX2}) because of possible discontinuities of
$Y^{\bolds{(\cdot)}}$. First notice that as $K\rightarrow \infty$, for $t$ fixed,
\[
\sum_{k=0}^{K-1}1_{(\tau_{k},\tau_{k+1}]}(t)
h_{\tau_k}\bigl(Y^{(\tau_k)},\overline{Z}_{\tau_k-}\bigr)
\lambda(\tau_k)\rightarrow
h_{t}\bigl(Y^{(t-)},\overline{Z}_{t-}\bigr) \lambda(t)
\]
for every $\omega\in\Omega$ fixed and for every $t< \tau$:
$Y^{\bolds{(\cdot)}}$ has limits from the left (Assumption~\ref{Ycadlag}), so
that as $\tau_k\uparrow t$, Restriction~\ref{hreg}(b) on $h$ can be used, and $\lambda$ is
continuous from the left [Assumption~\ref{int}(b)]. Taking integrals and
applying Lebesgue's dominated convergence theorem [$h_t$ and $\lambda$
are bounded because of Restriction~\ref{hreg}(a) and
Assumption~\ref{int}(a), resp.] leads to
\[
\sum_{k=0}^{K-1}
h_{\tau_k}\bigl(Y^{(\tau_k)},\overline{Z}_{\tau_k-}\bigr)
\lambda(\tau_k) (\tau_{k+1}-\tau_k)\rightarrow
\int_0^\tau  h_{t}\bigl(Y^{(t-)},\overline{Z}_{t-}\bigr) \lambda(t)\, dt
\]
for every $\omega\in\Omega$. As both $h$ and $\lambda$ are bounded,
Lebesgue's dominated convergence
theorem guarantees that indeed the expectation of the approximation
in (\ref{apY2}) converges to the expectation of the left-hand side of
(\ref{apY2}).
The same reasoning shows this for (\ref{apX2}), which
concludes the proof of equation~(\ref{seconde}) and
Theorem~\ref{thec}.  
\end{pf*}

\citet{Lok} shows that if the rest of the conditions in this section
are satisfied, Assumption~\ref{jump} (treatment does not
instantaneously affect the outcome of interest) is a necessary
condition for Theorem~\ref{thec}.
\begin{exa}[(Survival of AIDS patients and the Weibull proportional
hazards model)]\label{Wei}
\index{Weibull model}
Consider the setting of Examples~\ref{DAexa} and~\ref{DAexa2} and define
$N(t)=1$ if prophylaxis treatment started at or before time $t$ and
$0$ otherwise. Suppose that initiation of prophylaxis treatment can be
correctly modeled with the time-dependent Weibull proportional hazards
model
\[
\lambda_{\xi,\gamma,\theta}(t)=
1_{\{\mathrm{at}\ \mathrm{risk}\ \mathrm{at}\ t\}}\, \xi\gamma
t^{\gamma-1}e^{\theta_1 I_{\mathrm{AZT}}+\theta_2 I_{\mathrm{PCP}}(t)},
\]
where $I_{\mathrm{PCP}}(t)$ equals $1$ if the patient had PCP before time $t$ and
$0$ otherwise, and $\xi$ and $\gamma$ are greater than zero [for more
about the Weibull proportional hazards model and its applications see,
e.g., \citeauthor{Col} (\citeyear{Col})]. If the patient died before
$t$ or prophylaxis treatment already started before, the patient is
not ``at risk'' for initiation of treatment and, thus, $\lambda$ equals $0$.
Then the (partial) score equations for estimation of
$(\xi,\gamma,\theta)$ are
\[P_n
\int_0^\tau \pmatrix{
\dfrac{1}{\xi}&
\dfrac{1}{\gamma}+\log t&
I_{\mathrm{AZT}}&
I_{\mathrm{PCP}}(t)}^{\top}
\bigl(dN(t)- \lambda_{\xi,\gamma,\theta}(t)\,dt\bigr)=0.
\]
Such estimating equations can
also be written down for the model including $\alpha X_{\psi}$,
\[
\lambda_{\xi,\gamma,\theta,\alpha,\psi}(t)
=1_{\{\mathrm{at}\ \mathrm{risk}\ \mathrm{at}\ t\}}\,\xi\gamma
t^{\gamma-1}e^{\theta_1 I_{\mathrm{AZT}}+\theta_2 I_{\mathrm{PCP}}(t) + \alpha
X_{\psi}(t)}.
\]
\citet{Enc} proposes to estimate the
parameters in a model like this by choosing those parameters
$(\xi,\gamma,\theta,\psi)$ which maximize the likelihood
when $X_\psi$ is considered fixed and known, and for which
$\hat{\alpha}(\psi)=0$: for the true $\psi$, $X_\psi(t)=X(t)\sim Y^{(t)}$ does
not contribute to the model for treatment changes (under no unmeasured
confounding). To make the connection with the estimators in
the current article, notice that this leads to the same estimators as the
ones that solve the estimating equations arising from the likelihood
when $X_\psi$ is considered fixed and known, with $\alpha$ put to
zero. More precise, since we know that the true $\alpha$ is equal to $0$, we
put $\alpha$ equal to $0$ and get the estimating equations
\[P_n \int_0^\tau \pmatrix{\displaystyle
\frac{1}{\xi}&
\dfrac{1}{\gamma}+\log t&
I_{\mathrm{AZT}}&
I_{\mathrm{PCP}}(t)&
X_{\psi}(t)}^{\top}
\bigl(dN(t)-\lambda_{\xi,\gamma,\theta}(t)\,dt\bigr)=0\]
for the parameter $\psi$ (and thus also for $D$) and the
(nuisance) parameters $(\xi,\gamma,\theta)$.
These estimating equations are of the form of Theorem~\ref{thec},
\[
P_n\int_0^{\tau} h_t(X_\psi(t),\overline{Z}_{t-})
\bigl(dN(t)-\lambda_{\xi,\gamma,\theta}(t)\,dt\bigr)=0,
\]
but the function $h_t$ here is not bounded and $\lambda$ need not be
bounded (if $\gamma<1$), so unbiasedness does not follow immediately
from Theorem~\ref{thec}. However, we could restrict the interval
$[0,\tau]$ to $[\epsilon,\tau]$ for $\epsilon>0$
(to assure that $\lambda$ is bounded) and $\log t$ can be approximated
by the bounded functions $\log t\vee C$ ($C\rightarrow-\infty$) (to
make $h_t$ bounded), which all lead to unbiased estimating equations
because of Theorem~\ref{thec}. The above estimating equations are then
also unbiased because of Lebesgue's dominated convergence theorem [the
dominating function
is integrable since
\begin{eqnarray*}
E \int_0^\tau |\log t|
\bigl(dN(t)+\lambda_{\xi,\gamma,\theta}(t) \,dt\bigr)
&=&2 E \int_0^\tau |\log t| \lambda_{\xi,\gamma,\theta}(t)\, dt\\
&\leq & 2 \xi\gamma e^{|\theta_1|+|\theta_2|}
\int_0^\tau |\log t| t^{\gamma-1}\, dt,
\end{eqnarray*}
which is finite since $\gamma>0$].\vadjust{\goodbreak}

Under the
model for $D$ of Example~\ref{DAexa},
\begin{eqnarray*}
D_{\psi}(y,t;\overline{Z}_t)&=&(1-e^{\psi})1_{\{\mathrm{treated}\ \mathrm{at}\
t\}},\\
X_\psi(t)&=&t+\int_{t}^Y e^{\psi
1_{\{\mathrm{treated}\ \mathrm{at}\ s\}}}\, ds,\end{eqnarray*}
if the patient did not die before time $t$. In that case these are
five unbiased estimating equations for five unknown parameters. If the
parameter $\psi$ is of dimension greater than $1$, more unbiased
estimating equations can be constructed by adding more terms of the
form $\alpha f(X_{\psi}(t),\overline{Z}_{t-})$.
\end{exa}

\section{$X(t)$ does not predict treatment changes: a martingale result}
\label{mgs}

We show that, under no unmeasured confounding, just as $Y^{(t)}$,
$X(t)$ does not predict treatment changes, given past treatment- and
covariate history $\overline{Z}_{t-}$. We could hope for that since
$X(t)\sim Y^{(t)}$ given $\overline{Z}_t$ (see
Section~\ref{cmim}). The formal statement is (compare with
Assumption~\ref{intconf}, no unmeasured confounding) the following: the intensity
process $\lambda(t)$ of $N$ with respect to
$\sigma(\overline{Z}_t)$ is also the intensity process of
$N$ with respect to $\sigma(\overline{Z}_t,X(t))^a$. Then
$M(t)=N(t)-\int_0^t \lambda(s)\,ds$ is also a martingale with respect to
$\sigma(\overline{Z}_t,X(t))^a$. That will be useful later
when we study the behavior of estimators $\hat{\theta}$ and
$\hat{\psi}$ which are constructed with estimating equations of the
form of Theorem~\ref{thec},
$P_n\int_0^\tau h_t(X_{\psi}(t),\overline{Z}_{t-})
 (dN(t)-\lambda_{\theta}(t)\,dt)=0$.
For example, we can use the fact that usually
$\int_{[0,t]}H(s)\,dM(s)$ is a
martingale if $M$ is a martingale and $H$ a predictable process; a
sufficient condition for this is that $E\int
|H(s)|| dM(s)|<\infty$ [see,
e.g., Andersen et al. (\citeyear{ABGK})]. Hence, all estimating
equations of the above form which we saw before are in fact
martingales for
$(\theta,\psi)=(\theta_0,\psi_0)$.

Before going on, we first clarify why
$\sigma(\overline{Z}_t,X(t))$ is indeed a filtration. For
$s<t$, $X(s)$ is a deterministic (though unknown) function of
$(\overline{Z}_{t-},X(t))$ (i.e., if solutions to the
differential equation with $D$ are unique; see,
e.g., Theorem~\ref{isun} in the \hyperref[difeqap]{Appendix}). Similarly, for $s<t$, $X(t)$ is a
deterministic function of $(\overline{Z}_{t-},X(s))$. In
the rest of this article we will assume that these functions are
measurable functions on $\overline{\mathcal{Z}}_{t-}\times \mathbb{R}$
(sufficient conditions for that are that the infinitesimal
shift-function $D$ satisfies regularity Assumption~\ref{Dreg} below and that for
each $\omega\in\Omega$, $Z$ only jumps finitely many times; see
Appendix~\ref{smtb}, Lemma~\ref{mtb}). Thus,
%
\begin{equation}\label{fil}
\sigma(\overline{Z}_t,X(t))=\sigma\bigl(\overline{Z}_t,\bigl(X(s)\dvtx s\leq t\bigr)\bigr)
=\sigma(\overline{Z}_t,X(0)).
\end{equation}
We will use the filtration $\sigma(\overline{Z}_t,X(t))$ below,
keeping in mind that it is indeed a filtration and satisfies
equation~(\ref{fil}).

In the rest of this section we assume that the infinitesimal
shift-function $D$ satisfies the following regularity condition:
\begin{assu}[(Regularity of the infinitesimal
shift-function~$D$)] \label{Dreg}
\begin{longlist}[(a)]
\item[(a)] (\textit{Continuity between the jump times of} $Z$). If
$Z$ does not jump in $(t_1,t_2)$, then
$D(y,t;\overline{Z}_t)$ is continuous in $(y,t)$ on
$[t_1,t_2)$ and can be continuously extended to
$[t_1,t_2]$.
\item[(b)] (\textit{Boundedness}). For each $\omega\in\Omega$, there exists a constant
$C(\omega)$ such that
$|D(y,t;\overline{Z}_t)|\leq
C(\omega)$ for all $t\in[0,\tau]$ and all $y$.
\item[(c)] (\textit{Lipschitz continuity}). For
each $\omega\in\Omega$, there exist constants $L_1(\omega)$
and $L_2(\omega)$ with
\[
|D(y,t;\overline{Z}_t)-D(z,t;\overline{Z}_t)|\leq
L_1(\omega)|y-z|
\]
for all $t\in[0,\tau]$ and all $y,z$ and
\[
|D(y,t;\overline{Z}_t)-D(y,s;\overline{Z}_s)|\leq
L_2(\omega)|t-s|
\]
if $s<t$ and $Z$ does not jump in $(s,t]$.
\end{longlist}
\end{assu}

Most regularity conditions on $D$ here are satisfied for the $D$'s
from Appendix~\ref{XTap} [see also Lok (\citeyear{Lok,Mimarx})]. Only the second Lipschitz condition is extra. The
Lipschitz conditions are satisfied, for example, if, in between the jump
times of $Z$, $D$ is continuously differentiable with respect to $y$
and $t$ with derivatives which are bounded for every fixed
$\omega\in\Omega$.

The next theorem states that $M$ is indeed also a martingale
with respect to $\sigma(\overline{Z}_t,X(t))^a$:
\begin{theorem} \label{mgc}
Suppose that the conditions of Theorem~\textup{\ref{thec}} hold:
Assumptions~\textup{\ref{int}} (bounded intensity process), \textup{\ref{Ycadlag}} ($Y^{\bolds{(\cdot)}}$
cadlag), \textup{\ref{intconf}} (no unmeasured confounding), \textup{\ref{jump}} (no
instantaneous treatment effect) and for every
$t\in[0,\tau]$, $X(t)$ has the same distribution as $Y^{(t)}$
given $\overline{Z}_t$. Suppose, furthermore, that for each $\omega\in\Omega$, $Z$
jumps at most finitely many times, and that $D$ satisfies regularity
Condition~\ref{Dreg}.  Then the intensity process $\lambda(t)$ of $N$
with respect to $\sigma(\overline{Z}_t)$ is also the intensity process
of $N$ with respect to the filtration $\sigma(\overline{Z}_t,X(t))^a$.
\end{theorem}

Recall that in Section~\ref{cmim} we already mentioned that, under regularity
conditions, $X(t)$ mimics $Y^{(t)}$ in the sense that it has the same
distribution as $Y^{(t)}$ given $\overline{Z}_t$.

Under local rank preservation (see Section~\ref{loc}), $X(t)=Y^{(t)}$
for each $t$. In that case Theorem~\ref{mgc} would be the same as the
Assumption of no unmeasured confounding~\ref{intconf}. However, as
argued in Section~\ref{loc}, local rank preservation is generally
considered implausible.
\begin{pf*}{Proof of Theorem \ref{mgc}}  Because of Assumption~\ref{int},
$\Lambda(t)=\int_0^t \lambda(s)\, ds$ is predictable with respect to
$\sigma(\overline{Z}_t)$, so then it is also predictable
with respect to the larger filtration
$\sigma(\overline{Z}_t,X(t))^a$. We still need to prove that
$M$ is a martingale with respect to
$\sigma(\overline{Z}_t,X(t))^a$. Since a cadlag martingale
with respect to some filtration is also a martingale with respect to
its usual augmentation [see \citeauthor{RW} (\citeyear{RW}),
Lemma~67.10], it suffices to prove that $M$ is a martingale with
respect to $\sigma(\overline{Z}_t,X(t))$. Thus we need to
prove that, for $t_2>t_1$,
\[
E[M(t_2)-M(t_1)|\overline{Z}_{t_1},X(t_1)]
=0. \]
This is not immediate, since we do not have or expect that
$(X(t_1),\overline{Z}_{t_2})\sim
(Y^{(t_1)},\overline{Z}_{t_2})$ if $t_1<t_2$.

By the definition of conditional expectation, the above is the same as
\begin{equation}
\int_B \bigl(M(t_2)-M(t_1)\bigr)\,dP =0
\label{enul}
\end{equation}
for all $B\in\sigma(\overline{Z}_{t_1},X(t_1))$. Because
of Theorem~34.1 in \citet{Bill86}, it is sufficient to consider $B$'s
forming a $\pi$-system generating
$\sigma(\overline{Z}_{t_1},X(t_1))$. With $\sigma_1$ the
$\sigma$-algebra on $\overline{\mathcal{Z}}_{t_1}$,
\[
\{\omega\in\Omega\dvtx \overline{Z}_{t_1}\in A\mbox{ and }
X(t_1)\in(x_1,x_2)\dvtx A\in\sigma_1 \mbox{ and } x_1<x_2\in\mathbb{R}\}
\]
is such a $\pi$-system: it is closed under the formation of finite
intersections and generates
$\sigma(\overline{Z}_{t_1},X(t_1))$.  Therefore, we only
consider $B$'s of this form. We prove (\ref{enul}) for any
$B=\{\overline{Z}_{t_1}\in A\}\cap\{
X(t_1)\in(x_1,x_2)\}$. Let
$1^{(n)}_{(x_1,x_2)}$ be any approximation of
$1_{(x_1,x_2)}$ which is continuous for every fixed $n$,
with $1^{(n)}_{(x_1,x_2)}(x)\rightarrow
1_{(x_1,x_2)}(x)$ for every $x$ as $n\rightarrow \infty$ and
$|1^{(n)}_{(x_1,x_2)}|\leq 1$ for all $x$ and
$n$. Then
\begin{eqnarray*}
&&\int_B \bigl(M(t_2)-M(t_1)\bigr)\,dP\\
&&\qquad=E\bigl(1_B \cdot \bigl(M(t_2)-M(t_1)\bigr)\bigr)\\
&&\qquad=E \biggl(1_A(\overline{Z}_{t_1})
1_{(x_1,x_2)}(X(t_1))
\int_{(t_1,t_2]}\, dM(t)\biggr)\\
&&\qquad=E\int 1_{(t_1,t_2]}(t)\, 1_A(\overline{Z}_{t_1})
1_{(x_1,x_2)}(X(t_1))\, dM(t)\\
&&\qquad=E \int 1_{(t_1,t_2]}(t)\,
1_A(\overline{Z}_{t_1})
\lim_{n\rightarrow \infty}
1^{(n)}_{(x_1,x_2)}(X(t_1))\, dM(t)\\
&&\qquad=E\lim_{n\rightarrow \infty} \int 1_{(t_1,t_2]}(t)\,
1_A(\overline{Z}_{t_1})
1^{(n)}_{(x_1,x_2)}(X(t_1))
\bigl(dN(t)-\lambda(t)\,dt\bigr)\\
&&\qquad=\lim_{n\rightarrow \infty} E \int 1_{(t_1,t_2]}(t)\,
1_A(\overline{Z}_{t_1})
1^{(n)}_{(x_1,x_2)}(X(t_1))\,dM(t).
\end{eqnarray*}
The last two equalities follow from Lebesgue's dominated convergence
theorem [the prior to last equality since, for $\omega\in \Omega$
fixed, the integral is bounded since $N$ is finite and $\lambda$ is
bounded; the last equality since the integrals are all bounded
by $N(\tau)+\int_0^\tau \lambda(t)\,dt$, whose expectation is bounded by
$2\tau$ times the upper bound of~$\lambda$]. Equation~(\ref{enul}) and
the result of the theorem would follow from Theorem~\ref{thec} if
\[ 1_{(t_1,t_2]}(t)\, 1_A(\overline{Z}_{t_1})
1^{(n)}_{(x_1,x_2)}(X(t_1)) =
h^{(n)}_t(X(t),\overline{Z}_{t-})\]
for some $h^{(n)}_t$ from $\mathbb{R}\times \overline{\mathcal{Z}}_{t-}\rightarrow \mathbb{R}$
satisfying Restriction~\ref{hreg} for
each fixed $n$. In principle, this seems possible, since
$X(t_1)$ is a function of $X(t)$ and $\overline{Z}_{t-}$
for every $t>t_1$.

Indeed, under the conditions above on $D$ and $Z$, it is possible to
find such an $h^{(n)}_t$, as follows. Write
$x(\cdot;t_0,x_0)$ for the solution of the differential
equation
\[x{}'(t)=D(x(t),t;\overline{Z}_t)
\]
with (final or initial, depending on $t$) condition
$x(t_0)=x_0$. Existence and uniqueness of
$x(\cdot;t_0,x_0)$ on $[0,\tau]$ for every fixed
$\omega\in\Omega$ follows from Theorem~\ref{isun} in Appendix~\ref{difeqap}.
In this notation,
\begin{eqnarray*}
&&1_{(t_1,t_2]}(t)\,
1_A(\overline{Z}_{t_1})
1^{(n)}_{(x_1,x_2)}(X(t_1))\nonumber\\
&&\qquad=1_{(t_1,t_2]}(t)\,
1_A(\overline{Z}_{t_1})
1^{(n)}_{(x_1,x_2)}(x(t_1;t,X(t)))
= h^{(n)}_t(X(t),\overline{Z}_{t-})
\end{eqnarray*}
with
\begin{equation}\label{ht}
h^{(n)}_t(y,\overline{Z}_{t-})=1_{(t_1,t_2]}(t)\,
1_A(\overline{Z}_{t_1})
1^{(n)}_{(x_1,x_2)}(x(t_1;t,y)).
\end{equation}
We have to show that (\ref{ht}) satisfies
Restriction~\ref{hreg}. First we show that, for $t$ fixed,
$h^{(n)}_t:\mathbb{R}\times \overline{\mathcal{Z}}_{t-}\rightarrow
\mathbb{R}$ is measurable. From (\ref{ht}) we see that this is the
case if $x(t_1;t,\cdot)\dvtx\mathbb{R}\times \overline{\mathcal{Z}}_{t-}\rightarrow
\mathbb{R}$ is measurable, which follows
immediately from Lemma~\ref{mtb}. Restriction~\ref{hreg}(a) is
immediate, since $h^{(n)}_t$ is bounded by 1. For
Restriction~\ref{hreg}(b), we have to prove that, for all
$\omega\in\Omega$,
$h^{(n)}_t(y,\overline{Z}_{t-}(\omega))\rightarrow
h^{(n)}_{t_0}(y_0,\overline{Z}_{t_0-}(\omega))$ when
$y\rightarrow y_0$ and $t\uparrow t_0$. Fix $\omega\in\Omega$. We
consider three different kinds of $t_0$. If $t_0\leq t_1$ and
$t\uparrow t_0$, $h^{(n)}_t()=0=h^{(n)}_{t_0}(\cdot)$, so that the
convergence follows immediately.  If $t_0> t_2$ and $t\uparrow t_0$,
eventually $h^{(n)}_t(\cdot)=0=h^{(n)}_{t_0}(\cdot)$, so that the convergence
also follows immediately. If $t_0\in(t_1,t_2]$, convergence
of the first two factors is immediate.
For the last factor, we need differential equation theory.
$1_{(x_1,x_2)}^{(n)}$ is continuous. Thus, to prove that the
last factor in equation~(\ref{ht}) converges, it suffices to show that
$x(t_1;t,y)\rightarrow x(t_1;t_0,y_0)$ as
$t\uparrow t_0$ and $y\rightarrow y_0$.

Fix $\omega\in\Omega$. For $t$ close enough to $t_0$, we compare the
solution of the differential equation with final condition $y$ at $t$
with the solution of the differential equation with final condition
$y_0$ at $t_0$; we look at the value of the solution at the time point
$t_1$ before both $t$ and $t_0$. First, notice that because of
existence and uniqueness of solutions (Theorem~\ref{isun}), the
solution of the differential equation with final condition $y$ at $t$
takes a unique value $\tilde{y}=x(t_0;t,y)$ at $t_0$. Since
$x$ is differentiable with respect to its first argument with derivative $D$ and
$D$ is bounded by $C(\omega)$ [Assumption~\ref{Dreg}(b)],
$\tilde{y}$ is not far from $y$ if $t$ is not far from $t_0$:
\begin{equation}\label{yyt}
|\tilde{y}-y|=|x(t_0;t,y)-y|
\leq C(\omega) |t-t_0|.
\end{equation}
Next, notice that, again because of existence and uniqueness of
solutions, the value at $t_1$ of the solution of the differential
equation with final condition $y$ at $t$ is the same as the value at
$t_1$ of the solution of the differential equation with final
condition $\tilde{y}=x(t_0;t,y)$ at $t_0$.
This observation
implies that
\begin{eqnarray}\label{xafs}
|x(t_1;t,y)-x(t_1;t_0,y_0)|
&=&|x(t_1;t_0,\tilde{y})-x(t_1;t_0,y_0)|\nonumber\\
&\leq& e^{L_1(\omega)|t_1-t_0|}
|\tilde{y}-y_0|\nonumber\\[-8pt]\\[-8pt]
&\leq& e^{L_1(\omega)|t_1-t_0|}
(|\tilde{y}-y|+|y-y_0|)\nonumber\\
&\leq& e^{L_1(\omega)|t_1-t_0|}
\bigl(C(\omega) |t-t_0|+|y-y_0|\bigr).\nonumber
\end{eqnarray}
For the first inequality, we use Corollary~\ref{difbc1} and
Assumption~\ref{Dreg} [notice that possible jumps of $D$ at the
jump times of $Z$ do not matter here since one can split up the
interval, so if, e.g., there is just one jump at
$\tilde{t}\in(t_1,t_0]$, one gets a factor
\[
e^{L_1(\omega)|t_1-\tilde{t}|} \cdot
e^{L_1(\omega)|\tilde{t}-t_0|}=
e^{L_1(\omega)|t_1-t_0|},
\]
etc.\ (a formal proof can be given with induction since, with
$\omega\in\Omega$ still fixed, there are only finitely many jumps of
$Z$)]. For the last inequality, we use equation~(\ref{yyt}). If $y\rightarrow y_0$
and $t\uparrow t_0$, the bound in equation~(\ref{xafs}) converges to $0$ for
every fixed $\omega\in\Omega$.  Thus, indeed, if $y\rightarrow y_0$ and
$t\uparrow t_0$, $x(t_1;t,y)$ converges to
$x(t_1;t_0,y_0)$.  This finishes the proof.
\end{pf*}

\section{Consistency and asymptotic normality}
\label{asnorm}

The estimating equations for $(\theta,\psi)$ from
Section~\ref{ees} were all of the form $P_n
g_{\theta,\psi}(Y,\overline{Z})=0$. In the current section
we choose the dimension of $g$ the same as the dimension of
$(\theta,\psi)$. Estimating equations of this form are
well known. Theorem~\ref{consas} below is an example of asymptotic
theory in the setting of this article, with conditions in terms of $h$
and the intensity process $\lambda$.  Notice, however, that these
conditions are in fact stronger than necessary. For more theory about
these types of estimating equations and less restrictive conditions, see
\citet{Vaart}, Chapter~5.  In particular, conditions could be weakened
by considering the estimating equations as a whole instead of looking
at $h$ and $\lambda$ separately (see, e.g., Example~\ref{Weilem}).


We only consider smooth $h_{t}^{\theta,\psi}$:
\begin{rest} \label{hregpc} The functions
  $h_{t}^{\theta,\psi}\dvtx\mathbb{R}\times \overline{\mathcal{Z}}_{t}\rightarrow \mathbb{R}^m$ are measurable
  and:
\begin{longlist}[(a)]
\item[(a)] Every component of $h_{t}^{\theta_0,\psi_0}$ satisfies
Restriction~\ref{hreg}.
\item[(b)] $h_{t}^{\theta,\psi}(y,\overline{Z}_{t-})$ is bounded
  by a constant $C_1$ not depending on $\theta$, $\psi$, $t$, $y$ and
  $\omega\in\Omega$.
\item[(c)] For each $t\in[0,\tau]$ and $\omega\in\Omega$,
  $(\theta,\psi,y)\rightarrow
  h_{t}^{\theta,\psi}(y,\overline{Z}_{t-})$ is continuous.
\item[(d)] There exists a neighborhood of $(\theta_0,\psi_0)$
  such that, for each $t\in[0,\tau]$ and $\omega\in\Omega$,
  $h_{t}^{\theta,\psi}(y,\overline{Z}_{t-})$ is
  continuously differentiable with respect to $\theta$, $\psi$ and
  $y$, and these derivatives are all bounded by a constant $C_2$ not
  depending on $\theta$, $\psi$, $t$, $y$ and $\omega\in\Omega$.
\item[(e)] Every component of
  $\frac{\partial}{\partial (\theta,\psi)}|_{(\theta,\psi)=(\theta_0,\psi_0)}
  h_{t}^{\theta,\psi}(y,\overline{Z}_{t-})$ satisfies Restriction~\ref{hreg}.
\end{longlist}
\end{rest}
\begin{theorem}[(\textup{Consistency and asymptotic normality})]\label{consas}
   Suppose that
  Assumptions~\textup{\ref{int}} (bounded intensity process), \textup{\ref{Ycadlag}}
  [$Y^{(\bolds{\cdot})}$ cadlag], \textup{\ref{intconf}} (no unmeasured confounding)
  and~\textup{\ref{jump}} \textup{(}no instantaneous treatment effect\textup{)} are satisfied.
  Suppose also that, for every $t\in[0,\tau]$, $X(t)$ has
  the same distribution as $Y^{(t)}$ given $\overline{Z}_t$. From
  Theorem~\textup{\ref{thec}} we know that, for $h$ satisfying
  Restriction~\textup{\ref{hregpc}(a)}, $(\theta_0,\psi_0)$ is a zero
  of
\[E \int_0^\tau
h^{\theta,\psi}_t(X_\psi(t),\overline{Z}_{t-})
 \bigl(dN(t)-\lambda_\theta(t)\,dt\bigr).
\]
Suppose now that $(\theta_0,\psi_0)$ is the only zero.
Suppose, furthermore, that we know that
$(\theta_0,\psi_0)\in(\Theta,\Psi)$ with
$(\Theta,\Psi)$ compact, that $\theta\rightarrow
\lambda_\theta(t)$ is continuous for each $t$ and bounded
by a constant $C_3$ which does not depend on $(\omega,t,\theta)$,
and that $\psi\rightarrow X_\psi(t)$ is continuous for each $t$. Then
any sequence of (almost) zeros $(\hat{\theta},\hat{\psi})$
of
\[
\Psi_n(\theta,\psi)=
P_n\int_0^\tau h^{\theta,\psi}_t(X_\psi(t),\overline{Z}_{t-})
 \bigl(dN(t)-\lambda_\theta(t)\,dt\bigr),
 \]
that is, any sequence of estimators $(\hat{\theta},\hat{\psi})$ such that
$\Psi_n(\hat{\theta},\hat{\psi})$ converges in probability to zero, is a consistent
estimator for $(\theta_0,\psi_0)$ for each $h_t$ satisfying Restriction~\textup{\ref{hregpc}(a)--(c)}.

Suppose, moreover, that $\theta\rightarrow \lambda_\theta(t)$ is
differentiable with the respect to $\theta$ with derivative bounded by a
constant $C_4$ in a neighborhood of $\theta_0$, and $\psi\rightarrow
X_\psi(t)$ is differentiable with respect to $\psi$ with the derivative
bounded by a constant $C_5$ in a neighborhood of $\psi_0$.
Then for each $h$ satisfying Restriction~\ref{hregpc} there is a
neighborhood of $(\theta_0,\psi_0)$ such that
$E\int_0^\tau
h_t^{\theta,\psi}(X_\psi(t),\overline{Z}_{t-})
(dN(t)-\lambda_\theta(t)\,dt)$ is continuously differentiable
with respect to $(\theta,\psi)$. Suppose, moreover, that the
matrix
\[ V_0=E\biggl(\frac{\partial}{\partial (\theta,\psi)}\bigg|_{(\theta,\psi)=(\theta_0,\psi_0)}
\int_0^\tau
h^{\theta,\psi}_t(X_\psi(t),\overline{Z}_{t-})
 \bigl(dN(t)-\lambda_\theta(t)\,dt\bigr)\biggr)
\]
is nonsingular. Then there exists a sequence of (almost) zeros
$(\hat{\theta},\hat{\psi})$ to
$\Psi_n(\theta,\psi)$. Furthermore, any such sequence is
asymptotically normal:
\begin{equation}\sqrt{n}\bigl(\pmatrix{
\hat{\theta}&
\hat{\psi}}^\top
-\pmatrix{\theta_0&
\psi_0}^\top\bigr) \leadsto \mathcal{N}(0,V_0^{-1}W_0(V_0^{-1})^\top)
\label{asvar}\end{equation}
with $V_0$ the matrix above and, with $a^{\otimes 2}=a a^\top$,
\[W_0
=E\biggl(\biggl(\int_0^\tau
h^{\theta_0,\psi_0}_{t}(X_{\psi_0}(t),\overline{Z}_{t-})
\bigl(dN(t)-\lambda_{\theta_0}(t)\bigr)\biggr)^{\otimes 2}\biggr).
\]
\end{theorem}
\begin{pf*}{Proof} Consistency follows from Theorem~5.9 of \citet{Vaart}.
Repeatedly applying Lebesgue's dominated convergence theorem shows that our
conditions imply the conditions of the first and second paragraph after this Theorem~5.9.

The existence of (almost) zeros follows from \citet{VW},
Section~3.9,
Problem~9, whose solution is practically given by the hint below it.
This Problem~9 states that if $f\dvtx\Theta\times\Psi\rightarrow {\mathbb R}^d$
is a homeomorphism of a neighborhood of
$(\theta_0,\psi_0)\in\mathbb{R}^d$ onto a neighborhood of
$0\in {\mathbb R}^d$, then every continuous
$f\dvtx\Theta\times\Psi\rightarrow \mathbb{R}^d$ for which
$\sup_{(\theta,\psi)\in\Theta\times\Psi}\|
  f(\theta,\psi)-g(\theta,\psi)\|$ is
sufficiently small has at least one zero. In our case
$g(\theta,\psi)=E\int_0^\tau
h_t^{\theta,\psi}(X_\psi(t),\overline{Z}_{t-})
(dN(t)-\lambda_\theta(t))$ is continuously differentiable
in a neighborhood of $(\theta_0,\psi_0)$ by
Restriction~\ref{hregpc}(d) and the assumptions on $\lambda_\theta$ and
$X_\psi(t)$, under which differentiation and integration can be
exchanged (twice). The derivative of this $g(\theta,\psi)$
at $(\theta_0,\psi_0)$ is nonsingular by assumption, and
hence, $g$ is a homeomorphism of a neighborhood of
$(\theta_0,\psi_0)\in\mathbb{R}^d$ onto a neighborhood of
$0\in\mathbb{R}^d$. $\Psi_n(\theta,\psi)$ is continuous in
$(\theta,\psi)$ and close enough to
$g(\theta,\psi)$ for large $n$ with probability approaching
$1$ because of the second paragraph below Theorem~5.9 in \citet{Vaart}.
Hence, $\Psi_n(\theta,\psi)$ has a zero with probability
approaching $1$.

Asymptotic normality follows from Theorem~5.21 of \citet{Vaart}, as follows.
Define $g_{\theta,\psi}(Y,\overline{Z}_t)=\int_0^\tau h_t^{\theta,\psi}(X_\psi(t),
\overline{Z}_{t-})(dN(t)-\lambda_\theta(t)\,dt)$, which is continuously differentiable
with respect to $(\theta,\psi)$ in a neighborhood $U$ of $(\theta_0,\psi_0)$ under
our conditions. Making $U$ smaller so that all boundedness conditions hold on $U$,
we define $\dot{g}(Y,\overline{Z}_{t})=
\sup_{(\theta,\psi)\in U}\|\frac{\partial}{\partial (\theta,\psi)}g_{\theta,\psi}(Y,\overline{Z}_{t})\|$,
which is bounded by
$(C_2+C_2 C_5 + C_2) (N(\tau)+C_3\tau)+C_1 C_4 \tau$,
a constant plus a constant times $N(\tau)$.  This
$\dot{g}(Y,\overline{Z}_{t})$ is square integrable since
$N(\tau)$ is square integrable: it is well known that counting
processes with bounded intensity processes are square integrable [it
follows from, e.g., Proposition~II.4.1 of Andersen et al.
(\citeyear{ABGK})].
For the same reason,
$E\|g_{\theta,\psi}(Y,\overline{Z}_t)\|^2<\infty$.
The remaining conditions of Theorem~5.21 from \citet{Vaart} were
checked before, so, indeed, $(\hat{\theta},\hat{\psi})$ is
asymptotically normal with asymptotic covariance matrix (\ref{asvar}).
\end{pf*}

The asymptotic variance (\ref{asvar}) is often estimated by replacing
$(\theta_0,\psi_0)$ by their estimates and $E$ by $P_n$.
Thus, confidence intervals for $\psi_0$ can be constructed.  Also, tests
for whether $\psi_0$ has a specific value can be constructed that way.
For more about testing, see Section~\ref{test}.

One can often simplify the expression for the asymptotic variance
in equation~(\ref{asvar}) using Corollary~\ref{casvar} below.
We use the following lemma:
\begin{lem} \label{caslem}
Suppose that the conditions of Theorem~\textup{\ref{mgc}} hold:
Assumptions~\textup{\ref{int} (}bounded intensity process\textup{), \ref{Ycadlag}} [$Y^{\bolds{(\cdot)}}$ cadlag],
\textup{\ref{intconf} (}no unmeasured confounding\textup{), \ref{jump} (}no
instantaneous treatment effect\textup{)}, for every $t\in[0,\tau]$,
$X(t)$ has the same distribution as $Y^{(t)}$ given $\overline{Z}_t$ (see
Section~\textup{\ref{cmim}}), for each $\omega\in\Omega$, $Z$ jumps at most
finitely many times, and $D$ satisfies regularity
Condition~\textup{\ref{Dreg}}. Write $a\otimes b=a b^\top$ and
$a^{\otimes 2}=a a^\top$. Then for $h_t\dvtx\mathbb{R}\times
\overline{\mathcal{Z}}_{t}\rightarrow \mathbb{R}^m$ every component of
which satisfies Restriction~\textup{\ref{hreg},}
\[
E\biggl(\biggl(\int_0^\tau
h^{\theta_0,\psi_0}_{t}(X(t),\overline{Z}_{t-})\,dM(t)\biggr)^{\otimes 2}\biggr)
=E\int_0^\tau
h^{\theta_0,\psi_0}_{t}(X(t),\overline{Z}_{t-})^{\otimes 2}
\lambda(t)\,dt.
\]
If, furthermore, $\lambda_{\theta}$ is a correctly specified model for
$\lambda$ such that $\frac{\partial}{\partial \theta}\lambda_{\theta}$
exists and $D_{\psi}$ is a correctly specified model for $D$ such that, for
each $t$, $X_{\psi}(t)$ is differentiable with respect to $\psi$ at
$\psi=\psi_0$, then, for $h^{\theta,\psi}_{t}$ satisfying
Restriction~\ref{hregpc},
\begin{eqnarray*}
&&E\frac{\partial}{\partial
\theta}\bigg|_{\theta=\theta_0}
\int_0^\tau
h_{t}^{\theta,\psi_0}(X(t),\overline{Z}_{t-})
\bigl(dN(t)-\lambda_\theta(t)\,dt\bigr)\\
&&\qquad=-E\int_0^\tau
h^{\theta_0,\psi_0}_{t}
(X(t),\overline{Z}_{t-})
\cdot
\frac{\partial}{\partial \theta}\bigg|_{\theta=\theta_0}
\lambda_\theta(t)\,dt
\end{eqnarray*}
if the left- or right-hand side exists and
\begin{eqnarray*}
&&E\frac{\partial}{\partial
\psi}\bigg|_{\psi=\psi_0}
\int_0^\tau
h_{t}^{\theta_0,\psi}(X_{\psi}(t),\overline{Z}_{t-})\,dM(t)\\
&&\qquad=
E\int_0^\tau
\tilde{h}(X(t),\overline{Z}_{t-})
\otimes
\frac{\partial}{\partial \psi}\bigg|_{\psi=\psi_0}
X_{\psi}(t)\,dM(t)
\end{eqnarray*}
if the left- or right-hand side exists, where
$\tilde{h}(y,\overline{Z}_{t-})=\frac{\partial}{\partial
  y}h_t^{\theta_0,\psi_0}(y,\overline{Z}_{t-})$.
\end{lem}

\begin{pf*}{Proof}  For the first statement, we use counting process
theory from Andersen et al. (\citeyear{ABGK}), Chapter~2. If $M_1$ is a martingale,
$\langle M_1\rangle $ (if it exists) is defined as a predictable process
such that $M_1^2 - \langle M_1\rangle $ is a (local) martingale. If $M_2$
is another martingale, $\langle M_1,M_2\rangle $ (if it exists) is
defined as a predictable process such that $M_1 M_2
-\langle M_2,M_2\rangle $ is a (local) martingale. $\langle M_1\rangle $ is
called the predictable variation process of $M_1$ and
$\langle M_1,M_2\rangle $ is called the predictable covariation process
of $M_1$ and $M_2$. For vector-valued $M_1$, $\langle M_1\rangle $ is
defined as a predictable process such that $M_1^{\otimes
2}-\langle M_1\rangle $ is a (local) martingale. Hence, it is a matrix
with $\langle M_{1i},M_{1j}\rangle $ at the $i$th row, $j$th column.

As shown in Theorem~\ref{mgc}, $M(t)=N(t)-\int_0^t\lambda(s)\,ds$ is a
martingale with respect to the filtration
$\sigma(\overline{Z}_t,X(t))^a$. Counting process martingales like
this have compensators:
$\langle M(t)\rangle =\int_0^t \lambda(s)\,ds$ [see,
e.g., Proposition~II.4.1 in Andersen et al.
(\citeyear{ABGK})]. Moreover, if $H_1$ and $H_2$ are (locally) bounded
$\sigma(\overline{Z}_t,X(t))^a$-predictable processes with values in
$\mathbb{R}$, then
$\langle \int_0^t H_1(s)\, dM(s),\break \int_0^t H_2(s)\, dM(s)\rangle $ exists and
\[
\biggl<\int_0^t H_1(s)\, dM(s),\int_0^t H_2(s)\,
dM(s)\biggr>=\int_0^t H_1(s) H_2(s) \lambda(s)\,ds
\]
[Proposition~II.4.1 or (2.4.9) in Andersen et al.
(\citeyear{ABGK})]. Because $h_{t}$ satisfies Restriction~\ref{hreg},
$h_{t} (X(t),\overline{Z}_{t-})$ is a bounded
$\sigma(\overline{Z}_t,X(t))^a$-predictable process (proof just as in
Lemma~\ref{hpred}). Therefore, the theory above leads to
\begin{eqnarray*}
E\biggl(\biggl(\int_0^\tau
h^{\theta_0,\psi_0}_{t}(X(t),\overline{Z}_{t-})\,
dM(t)\biggr)^{\otimes 2}\biggr)
&=&E\biggl(\bigg<\int_0^\tau
h^{\theta_0,\psi_0}_{t}(X(t),\overline{Z}_{t-})\,
dM(t)\bigg>\biggr)\\
&=&E\biggl(\int_0^\tau
h^{\theta_0,\psi_0}_{t}(X(t),\overline{Z}_{t-})^{\otimes
2} \lambda(t)\,dt\biggr).
\end{eqnarray*}

For the second statement, notice that, under the conditions of the
lemma,
\begin{eqnarray*}
&&\frac{\partial}{\partial
\theta}\bigg|_{\theta=\theta_0}
\biggl(\int_{0}^{\tau}
h_{t}^{\theta,\psi_0}
(X(t),\overline{Z}_{t-})
\bigl(dN(t)-\lambda_\theta(t)\,dt \bigr)\biggr)\\
&&\qquad=-\int_{0}^{\tau}
h^{\theta_0,\psi_0}_{t}(X(t),\overline{Z}_{t-})
\otimes
\frac{\partial}{\partial
\theta}\bigg|_{\theta=\theta_0}\lambda_{\theta}(t)\,dt\\
&&\qquad\hspace*{18pt}{}+\int_{0}^{\tau}
\biggl(\frac{\partial}{\partial \theta}\bigg|_{\theta=\theta_0}
h_{t}^{\theta,\psi_0}\biggr)
(X(t),\overline{Z}_{t-})
\bigl(dN(t)-\lambda(t)\,dt\bigr),
\end{eqnarray*}
and the expectation of the second term here is equal to zero because
of Theorem~\ref{thec}.

For the third statement, notice that, under the conditions of the
lemma,
\begin{eqnarray*}
&&\frac{\partial}{\partial
\psi}\bigg|_{\psi=\psi_0}
\biggl(\int_{0}^{\tau}
h_{t}^{\theta_0,\psi}
(X_{\psi}(t),\overline{Z}_{t-})
\bigl(dN(t)-\lambda_\theta(t)\,dt \bigr)\biggr)\\
&&\qquad=
\int_{0}^{\tau}\biggl(
\tilde{h}_{t}(X(t),\overline{Z}_{t-})
\otimes
\frac{\partial}{\partial \psi}\bigg|_{\psi=\psi_0}
X_{\psi}(t)\biggr)
\bigl(dN(t)-\lambda(t)\,dt\bigr)\\
&&\qquad\hspace*{9pt}{}+\int_{0}^{\tau}
\biggl(\frac{\partial}{\partial \psi}\bigg|_{\psi=\psi_0}
h_{t}^{\theta_0,\psi}\biggr)
(X(t),\overline{Z}_{t-})
\bigl(dN(t)-\lambda(t)\,dt \bigr)
\end{eqnarray*}
because of the chain rule, and the expectation of
the second term is equal to zero because of Theorem~\ref{thec}.
\end{pf*}

 This lemma simplifies the asymptotic
variance formula of the estimators in equation~(\ref{asvar}):
\begin{cor}[(\textup{Asymptotic variance})]
\label{casvar}
Suppose that the conditions of Theorem~\textup{\ref{mgc}} hold:
Assumptions~\textup{\ref{int}} (bounded intensity process), \textup{\ref{Ycadlag}} [$Y^{\bolds{(\cdot)}}$
cadlag], \textup{\ref{intconf} (}no unmeasured confounding\textup{), \ref{jump}} (no
instantaneous treatment effect), for every $t\in[0,\tau]$,
$X(t)$ has the same distribution as $Y^{(t)}$ given $\overline{Z}_t$ (see
Section~\textup{\ref{cmim}}), for each $\omega\in\Omega$, $Z$ jumps at most
finitely many times, and $D$ satisfies regularity
Condition~\textup{\ref{Dreg}}. Suppose also that $\lambda_{\theta}$ is a
correctly specified model for $\lambda$ such that
$\frac{\partial}{\partial \theta}\lambda_{\theta}$ exists and that
$D_{\psi}$ is a correctly specified model for $D$ such that, for each
$t$, $X_{\psi}(t)$ is differentiable with respect to $\psi$ at
$\psi=\psi_0$.  Then if
\[g_{\theta,\psi}(Y,\overline{Z})=\int_0^\tau
h^{\theta,\psi}_{t}(X_{\psi}(t),\overline{Z}_{t-})
\bigl(dN(t)-\lambda_\theta(t)\,dt\bigr)
\]
and $h^{\theta,\psi}_{t}$ satisfies Restriction~\ref{hregpc}, the
asymptotic variance (\ref{asvar}) is equal to
$V_0^{-1}W_0 {V_0^{-1}}^\top$ with
\[
W_0=E\biggl(\int_0^\tau h_{t}^{\theta_0,\psi_0}
(X(t),\overline{Z}_{t-})^{\otimes 2}
\lambda(t)\,dt\biggr)
\]
and $V_0=(V_{0\theta}V_{0\psi})$ with
\[
V_{0\theta}=-E\biggl(\int_0^\tau\biggl(
h^{\theta_0,\psi_0}_{t}
(X(t),\overline{Z}_{t-})
\otimes
\frac{\partial}{\partial \theta}\bigg|_{\theta=\theta_0}
\lambda_\theta(t)\biggr)\,dt\biggr)
\]
and, with
$\tilde{h}(y,\overline{Z}_{t-})=\frac{\partial}{\partial
  y}h_t^{\theta_0,\psi_0}(y,\overline{Z}_{t-})$,
\[
V_{0\psi}=E\biggl(\int_0^\tau\biggl(
\tilde{h}(X(t),\overline{Z}_{t-})
\otimes
\frac{\partial}{\partial \psi}\bigg|_{\psi=\psi_0}
X_{\psi}(t)\biggr)
\bigl(dN(t)-\lambda(t)\,dt
\bigr)\biggr).
\]
\end{cor}

We conclude this section with an example to see the machinery work in
practice. Notice that the boundedness conditions of
Theorem~\ref{consas} are somewhat too restrictive for this example,
but that the results hold true under these weaker restrictions, too.
\begin{lem}
[\textup{(Survival of AIDS patients and the Weibull proportional hazards
model)}] \label{Weilem}  \index{Weibull model} Consider the setting of
Example~\textup{\ref{Wei}}, and suppose that the assumptions of
Section~\ref{mgs} are satisfied. In Example~\textup{\ref{Wei},}
\[
\lambda_{\xi,\gamma,\theta}(t)=
1_{\{\mathrm{at}\ \mathrm{risk}\ \mathrm{at}\ t\}}\xi\gamma
t^{\gamma-1}e^{\theta_1 I_{\mathrm{AZT}}+\theta_2 I_{\mathrm{PCP}}(t)},
\] where
at risk means at risk for initiation of prophylaxis treatment.
$X_{\psi}(t)$ is the solution to the differential equation
$X_{\psi}'(t)=D_{\psi}(X_{\psi}(t),t;\overline{Z}_t)$ with final
condition $X_{\psi}(\tau)=Y$ and
$D_{\psi}(y,t;\overline{Z}_t)=(1-e^\psi)
1_{\{\mathrm{treated}\mathrm{at}t\}}$, so
\[
X_{\psi}(t)=t
+\int_{t}^Y e^{\psi1_{\{\mathrm{treated}\ \mathrm{at}\ s\}}}\,ds.
\]
In Example~\ref{Wei} we already saw that
$(\xi_0,\gamma_0,\theta_0,\psi_0)$ is a zero of
\[
E\biggl( \int_0^\tau
\pmatrix{\dfrac{1}{\xi}&
\dfrac{1}{\gamma}+\log t&
I_{\mathrm{AZT}}&
I_{\mathrm{PCP}}(t)&
X_{\psi}(t)}^{\top}
\bigl(dN(t)-\lambda_{\xi,\gamma,\theta}(t)\,dt\bigr)\biggr).
\]
Suppose now that $(\xi_0,\gamma_0,\theta_0,\psi_0)$ is the
only zero.
Suppose, furthermore, that the survival time $Y$ takes values in a
compact space $[0,y_0]\subset\mathbb{R}$ and that we know
that $(\xi_0,\gamma_0,\theta_0,\psi_0)\in\Xi\times\Gamma
\times\Theta\times\Psi$, with $\Xi\subset(0,\infty)$,
$\Gamma\subset(0,\infty)$, $\Theta$ and $\Psi$ all four compact
(note that this implies that $\xi$ and $\gamma$ are bounded away
from $0$). Then any sequence of (almost) zeros
$(\hat{\xi},\hat{\gamma},\hat{\theta},\hat{\psi})$ of
\begin{eqnarray} \label{weiboh}
&&\Psi_n(\xi,\gamma,\theta,\psi)\nonumber\\
&&\qquad=P_n\int_0^\tau\pmatrix{
\dfrac{1}{\xi}&
\dfrac{1}{\gamma}+\log t&
I_{\mathrm{AZT}}&
I_{\mathrm{PCP}}(t)&
X_{\psi}(t)}^{\top}\\
&&\qquad \qquad\quad {}\times\bigl(dN(t)-\lambda_{\xi,\gamma,\theta}(t)\,dt\bigr),\nonumber
\end{eqnarray}
that is, any sequence of estimators
$(\hat{\xi},\hat{\gamma},\hat{\theta},\hat{\psi})$ such that
$\Psi_n(\hat{\xi},\hat{\gamma},\hat{\theta},\hat{\psi})$
converges in probability to zero, is a consistent estimator for
$(\xi_0,\gamma_0,\theta_0,\psi_0)$. Moreover,
$V_0=(V_{0\theta}V_{0\psi})$ as in
Corollary~\ref{casvar} exists, and
\begin{eqnarray*}
&&V_{0\theta}=-E\int_0^\tau
\pmatrix{\dfrac{1}{\xi_0}& \dfrac{1}{\gamma_0}+\log t& I_{\mathrm{AZT}}&
I_{\mathrm{PCP}}(t)&X(t)}^{\top}\\
&&\hspace*{60pt}{}\times\pmatrix{\dfrac{1}{\xi_0}& \dfrac{1}{\gamma_0}+\log t &I_{\mathrm{AZT}}&
I_{\mathrm{PCP}}(t)}\lambda(t)\, dt
\end{eqnarray*}
and $V_{0\psi}$ is a five-dimensional vector with zeros in the first
four positions and
\[
E\int_0^\tau \biggl(\frac{\partial}{\partial
\psi}\bigg|_{\psi=\psi_0} X_{\psi}(t)\biggr)
\bigl(dN(t)-\lambda(t)\,dt\bigr)
\]
in the fifth, with
\[
\frac{\partial}{\partial
\psi}\bigg|_{\psi=\psi_0}X_\psi(t)=\int_{t\wedge Y}^Y
e^{\psi_0}\,1_{\{\mathrm{treated}\ \mathrm{at}\ s\}}\, ds.
\]
If this $V_0$ is a nonsingular matrix,
then there exists a sequence of (almost) zeros
$(\hat{\xi},\hat{\gamma},\hat{\theta},\hat{\psi})$ of
(\ref{weiboh}). Furthermore, any such sequence is
asymptotically normal:
\[\sqrt{n}\bigl(\pmatrix{\hat{\xi}&
\hat{\gamma}&\hat{\theta}&
\hat{\psi}}^{\top}
-\pmatrix{\xi_0&\gamma_0& \theta_0 &\psi_0}^\top\bigr) \leadsto
\mathcal{N}\bigl(0,V_0^{-1}W_0(V_0^{-1})^\top\bigr),
\]
with $V_0$ the matrix above and
\[
W_0=
E\int_0^\tau\pmatrix{
\dfrac{1}{\xi_0}&
\dfrac{1}{\gamma_0}+\log t&
I_{\mathrm{AZT}}&
I_{\mathrm{PCP}}(t)&
X(t)}^{\top^{\otimes 2}}
\lambda(t)\, dt.
\]
Moreover, $\hat{\theta}$ and $\hat{\psi}$ are asymptotically
independent.
\end{lem}

The asymptotic independence here turns out to be no coincidence; see
Lok (\citeyear{Lok,plss}).
\begin{pf*}{Proof of Lemma \ref{Weilem}} The findings here are similar to the findings in
Theorem~\ref{consas}, but the boundedness conditions fail to hold here.
Consistency follows from \citet{Vaart}, Theorem~5.9. Existence of a
sequence of (almost) zeros follows from \citet{VW}, Section~3.9,
Problem~9, whose solution is practically given by the hint below it.
Asymptotic normality follows from \citet{Vaart}, Theorem~5.21.  The
asymptotic variance equals $V_0^{-1}W_0 (V_0^{-1})^\top$
because of Corollary~\ref{casvar} ($\lambda$ and $h$ are not bounded
here, but we can restrict the interval to $[\epsilon,\tau]$
for $\epsilon>0$ and let $\epsilon\downarrow 0$). We leave checking
the conditions of these theorems to the reader [or see
\citeauthor{Lok} (\citeyear{Lok}), Section~7.7].

Asymptotic independence of $\hat{\theta}$ and $\hat{\psi}$ follows by
direct calculation, after noticing that
\[
V_0=\pmatrix{-B&\mathbf{0}\cr
-C&A_{\psi_0}}
\Longrightarrow V_0^{-1}=\pmatrix{
-B^{-1} & \mathbf{0}\cr
- A_{\psi_0}^{-1}C
B^{-1} &
A_{\psi_0}^{-1}},
\]
with
\begin{eqnarray*}
A_{\psi_0}&=&E\int_0^\tau \biggl(\frac{\partial}{\partial
\psi}\bigg|_{\psi=\psi_0} X_{\psi}(t)\biggr)
\bigl(dN(t)-\lambda(t)\, dt\bigr),\\
B&=&E\int_0^\tau
\pmatrix{\dfrac{1}{\xi_0} &\dfrac{1}{\gamma_0}+\log t& I_{\mathrm{AZT}}& I_{\mathrm{PCP}}(t)}^{\top^{\otimes 2}}
\lambda(t)\, dt
\end{eqnarray*}
and
\[
C=E\int_0^\tau
X(t) \cdot
\pmatrix{\dfrac{1}{\xi_0}& \dfrac{1}{\gamma_0}+\log t &I_{\mathrm{AZT}} &I_{\mathrm{PCP}}(t)}
\lambda(t)\, dt.
\]
Define
\[
D=E\int_0^\tau X(t)^2\lambda(t)\,dt.
\]
This concludes the proof:
\begin{eqnarray*}
V_0^{-1}W_0(V_0^{-1})^\top
&=&\pmatrix{
-B^{-1} & 0\cr
-A_{\psi_0}^{-1}CB^{-1} & A_{\psi_0}^{-1}
}
\pmatrix{B&C^\top\cr
C&D}\\
&&{}\times \pmatrix{-B^{{-1}^\top} & -B^{{-1}^\top}C^\top A_{\psi_0}^{{-1}\top}\cr
0 & A_{\psi_0}^{{-1}^\top}}\\
&=&\pmatrix{(B^{-1})^{\top} & 0\cr
0 & A_{\psi_0}(CB^{-1}C^\top
+D)(A_{\psi_0}^{-1})^\top }.\end{eqnarray*}  
\end{pf*}

\section{Test for treatment effect without specifying a model for $D$}
\label{test}

We show that one can often test whether treatment affects the outcome
of interest without specifying a model $D_{\psi}$ for $D$. This was
conjectured, but not proved, in \citet{Enc}. If one does not have to
specify a model for $D$ in order to test whether treatment affects the
outcome, false conclusions caused by misspecification of the model for
$D$ can be avoided.


%

Under the null hypothesis of the no treatment effect, $D\equiv 0$ and
$X(t)\equiv Y$ (see Section~\ref{model}). If there is no unmeasured
confounding, $X(t)$ does not predict $N(t)$ given $\overline{Z}_{t-}$
(see Theorem~\ref{mgc}).  Hence, if there is no unmeasured confounding
and no treatment effect, adding the \textit{observed outcome} $Y$ to the
prediction model for treatment effect should not help the prediction.
This idea, presented in \citet{Enc} for the case of local rank
preservation (see Section~\ref{loc}), can be proven to be correct as
follows.

Technically, the tests in this section are similar to the score test
[for more about the score test see, e.g., \citeauthor{CH}
(\citeyear{CH})]. Suppose that the conditions of Section~\ref{mgs} are
satisfied, and that we have a correctly specified parametric model
$\lambda_{\theta}$ for $\lambda$. Define
\[
g_\theta(Y,\overline{Z})=\int_0^\tau
h_t^\theta(Y,\overline{Z}_{t-})
\bigl(dN(t)-\lambda_\theta(t)\,dt\bigr),
\]
with $h_t^{\theta_0}$ satisfying the regularity condition
Restriction~\ref{hreg}. The key idea of this procedure is that if
treatment does not affect the outcome, $D\equiv 0$, so $X(t)\equiv Y$,
and $g_{\theta_0}(Y,\overline{Z})$ has expectation zero
because of Theorem~\ref{thec}. Since $\theta_0$ is unknown, we base
the test on the limiting behavior under $D\equiv 0$ of $\sqrt{n}P_n
g_{\hat{\theta}}(Y,\overline{Z})$, where $\hat{\theta}$ is
an estimator of the nuisance parameter $\theta_0$. We will show that
if $D\equiv 0$, $\sqrt{n}P_n
g_{\hat{\theta}}(Y,\overline{Z})$ converges to a normal
random variable with expectation zero, which leads to a test for
whether $D\equiv 0$ in the usual way.

The nuisance parameter $\theta_0$ will be estimated using some set of
estimating equations
\[ P_n
\tilde{g}_\theta(\overline{Z})=0
\] with
$E\tilde{g}_{\theta_0}(\overline{Z})=0$,
$E\tilde{g}_{\theta}(\overline{Z})$ differentiable in
$\theta$ and
$E\tilde{g}_{\theta_0}^2(\overline{Z})<\infty$. A natural
choice would be a maximum (partial) likelihood estimator for
$\theta_0$. We suppose throughout this section that the resulting
estimator $\hat{\theta}$ is consistent and asymptotically normal with
\begin{equation}\label{thetaas}
\sqrt{n}(\hat{\theta}-\theta_0)=
-\biggl(\frac{\partial}{\partial \theta}\bigg|_{\theta=\theta_0}
E \tilde{g}_{\theta}(\overline{Z})\biggr)^{-1}
\sqrt{n}P_n\tilde{g}_{\theta_0}(\overline{Z}) + o_P(1),
\end{equation}
as will usually follow from, for example, \citet{Vaart}, Theorem~5.21.

If $\lambda_\theta$ and $h_t^\theta$ are sufficiently smooth,
\[
\theta\rightarrow \sqrt{n} P_n g_\theta(Y,\overline{Z})
\]
is differentiable with respect to $\theta$, and a Taylor expansion
around $\theta_0$ leads to
\[
\sqrt{n} P_n g_{\hat{\theta}}(Y,\overline{Z})
= \sqrt{n} P_n g_{\theta_0}(Y,\overline{Z})
+ P_n \dot{g}_{\tilde{\theta}}(Y,\overline{Z})
\sqrt{n}(\hat{\theta}-\theta_0),
\]
with $\dot{g}_{\theta}$ the derivative of $g_\theta$ with respect to
$\theta$
and $\tilde{\theta}$ between $\theta_0$ and $\hat{\theta}$. Since
$\hat{\theta}$ converges in probability to $\theta_0$, so does
$\tilde{\theta}$. Therefore, usually,
$
P_n \dot{g}_{\tilde{\theta}}(Y,\overline{Z}) \mathop{\stackrel{P}{\rightarrow}}\break E
\dot{g}_{\theta_0}(Y,\overline{Z})$.
Sufficient conditions under which this holds are given in
Appendix~\ref{GCap}, Lemma~\ref{pGC}.
Because of~(\ref{thetaas}) and the central limit theorem,
$\sqrt{n}(\hat{\theta}-\theta_0)$ converges in
distribution. Therefore, an application of Slutzky's lemma leads to
\begin{eqnarray*}
\sqrt{n}P_n g_{\hat{\theta}}(Y,\overline{Z})
&=&\sqrt{n} P_n g_{\theta_0}(Y,\overline{Z})
+E\dot{g}_{\theta_0}(Y,\overline{Z})\sqrt{n}(\hat{\theta}-\theta_0) +o_P(1)\\
&=&\sqrt{n} P_n g_{\theta_0}(Y,\overline{Z})
-E\dot{g}_{\theta_0}(Y,\overline{Z})V_0^{-1}\sqrt{n} P_n\tilde{g}_{\theta_0}(\overline{Z})
+o_P(1)\\
&=&(-E\dot{g}_{\theta_0}(Y,\overline{Z})
V_0^{-1}\ \
\mathrm{Id}_{\dim g_{\theta}})
\pmatrix{\sqrt{n}P_n\tilde{g}_{\theta_0}(\overline{Z})\cr
\sqrt{n}P_n g_{\theta_0}(Y,\overline{Z})}
+o_P(1),
\end{eqnarray*}
with $V_0=\frac{\partial}{\partial
\theta}|_{\theta=\theta_0}E\tilde{g}_\theta(\overline{Z})$. If $D\equiv 0$,
$X(t)\equiv Y$, so that Theorem~\ref{thec}
implies that also the expectation of
$g_{\theta_0}$ is equal to zero. Therefore, the central limit theorem can be
applied on the vector with $\sqrt{n}$ on the right-hand side;
it converges to a normal random variable with expectation
zero. Because of the Continuous Mapping Theorem [see, e.g., \citeauthor{Vaart}
(\citeyear{Vaart}), Chapter~18], $\sqrt{n}P_n
g_{\hat{\theta}}(Y,\overline{Z})$ then converges to a normal random
variable with expectation zero, too.

Calculation of its limiting covariance matrix is standard [see,
e.g.,\break
\citeauthor{Vaart} (\citeyear{Vaart}), Chapter~18]. To save space, we
omit that calculation here. If desirable, one can use
Theorem~\ref{thec} and Lemma~\ref{caslem} to simplify the expression.

Notice that a test for whether $D=D_0$ for any specific $D_0$ can be
constructed in exactly the same way. If we have a correctly
specified model $D_{\psi}$ for $D$, this thus also leads to a
confidence region for $\psi_0$ in the usual way, using the duality
between testing and confidence regions: include those $\psi$ for which
the null hypothesis $D=D_{\psi}$ is not rejected.

\section{Discussion and extensions}
\label{discu}

The proof of consistency and asymptotic normality of the estimators
presented in this article applies to continuous-time structural nested
models. A similar proof applies to structural nested models in
discrete time (when covariates are only measured at finitely many
fixed times $0=\tau_0<\tau_1<\cdots<\tau_K<\tau_{K+1}=\tau$, which are
the same for all patients and known in advance). \citet{SNart} argue
without proof that consistency and asymptotic normality should hold
for discrete time structural nested models under reasonable
assumptions; the proof is completed with the current article, as
follows. It is easy to see that in discrete time, $\sum_{\tau_k\leq
  t}P(\Delta
  N(\tau_k)=1|\overline{Z}_{\tau_k-})$ is the
compensator of $N$ with respect to $\sigma(\overline{Z}_t)$
[see, e.g., \citeauthor{Lok} (\citeyear{Lok}), Section~7.4]. The
assumption of no unmeasured confounding can be formalized as
\[P\bigl(\Delta
  N(\tau_k)=1|\overline{Z}_{\tau_k-},Y^{(\tau_k)}\bigr)=
P\bigl(\Delta N(\tau_k)=1|\overline{Z}_{\tau_k-}\bigr).
\]
Since $X(\tau_k)\sim Y^{(\tau_k)}$ given
$\overline{Z}_{\tau_k}$, this implies that also
\[P\bigl(\Delta
  N(\tau_k)=1|\overline{Z}_{\tau_k-},X(\tau_k)\bigr)=
P\bigl(\Delta N(\tau_k)=1|\overline{Z}_{\tau_k-}\bigr).
\]
The discrete-time counterparts of Theorem~\ref{thec} and
Theorem~\ref{mgc} follow immediately [see \citeauthor{Lok}
(\citeyear{Lok})]. Consistency and asymptotic normality follow in the
same way as for continuous time models.

The tests for treatment effect in this article can be carried out
without specifying a model for treatment effect; that is, no model for
$D$ is needed. This is an important feature of the tests because it
allows one to avoid false conclusions caused by misspecification of
the model for $D$. In practice, it may be hard to specify a correct
parametric model for the infinitesimal shift-function, $D$. Thus, it
is good that specification of a model for $D$ is not needed to test
for treatment effect.

The estimators in this article require the correct specification of a
model for treatment effect and of a model for prediction of treatment
changes. For the discrete-time setting, \citet{doub} has recently
proposed estimators which are doubly robust. Doubly robust estimators
are consistent and asymptotically normal if (i) the model for
prediction of treatment changes ($\lambda$ in the current article) is
correctly specified or if (ii) a regression model of a blipped down
outcome [$X_{\psi}(t)$ in the current article] on past treatment- and
covariate history $\overline{Z}_{t-}$ is correctly specified. In any
case, the model for treatment effect ($D_\psi$ in the current article)
has to be well specified.

In this article estimation started with the specification of a model
for the infinitesimal shift-function, $D$.  Interpretation of results
may be easier when one starts with a model like~(\ref{DAsim}),
$Y^{(t)}-t \sim \int_t^Ye^{\psi 1_{\{\mathrm{treated}\ \mathrm{at}\ s\}}}ds=e^{\psi}\cdot \mathit{DUR}(t,Y)+1\cdot
(Y-t-\mathit{DUR}(t,Y))$, given $\overline{Z}_{t}$. Here, $\mathit{DUR}(t,u)$ is the
duration of treatment in the interval $(t,u)$.  The main results in
the current article apply also to this model. The proofs for
Theorems~\ref{thec} and~\ref{consas} do not depend on $X$ being the
solution to $X'(t)=D(X(t),t;\overline{Z}_{t})$. The proof of
Theorem~\ref{mgc} does depend on $X'(t)=D(X(t),t;\overline{Z}_{t})$,
but it simplifies considerably if (\ref{DAsim}) or~(\ref{DAsim2}) is
used as a starting point.  Let us show this for (\ref{DAsim}). Define
$\tilde{X}(t)=t+\int_t^Ye^{\psi 1_{\{\mathrm{treated}\ \mathrm{at}\ s\}}}\,ds$. Using the first part of the proof of
Theorem~\ref{mgc}, for $t_1<t$,
\begin{eqnarray*}
&&1_{(t_1,t_2]}(t) 1_{A}(\overline{Z}_{t_1})1^{(n)}_{(x_1,x_2)}(\tilde{X}(t_1))\\
&&\qquad=1_{(t_1,t_2]}(t) 1_{A}(\overline{Z}_{t_1})1^{(n)}_{(x_1,x_2)}\bigl(\tilde{X}(t)+\mathit{DUR}(t_1,t)
(-1+e^\psi)\bigr)\\
&&\qquad=h_t^{(n)}(\tilde{X}(t),\overline{Z}_{t-}),
\end{eqnarray*}
where $\mathit{DUR}(t_1,t)$ is the duration of treatment in the interval
$(t_1,t)$. $h_t^{(n)}$ is measurable if the duration of
treatment until $t$ is in included in $Z(t)$. Moreover, if $t\uparrow
t_0$ and $y\rightarrow y_0$, then $h_t^{(n)}(y,t)\rightarrow
h_t^{(n)}(y_0,t_0)$. Hence, it follows immediately that
$h_t^{(n)}$ satisfies Restriction~\ref{hreg}, which concludes the
proof.
Depending on the specific application, it will be more appropriate to
start with $Y^{(t)}-t \sim \int_t^Ye^{\psi 1_{\{\mathrm{treated}\ \mathrm{at}\ s\}}}\,ds$ given $\overline{Z}_{t}$ or
with a model for $D$; see, for example, Example~\ref{prior}.

In the previous literature [see, e.g., Robins et al. (\citeyear{Aids}), \citet{smoke},
Witteman et al. (\citeyear{hyp}), Keiding et al. (\citeyear{Keiding}) and Hern\'{a}n et al.
(\citeyear{Her})]
 applications have
been carried out under the assumption of local rank preservation,
where $X(t)=Y^{(t)}$ for each $t$. As pointed out in these papers, and
as discussed in Section~\ref{loc}, the assumption of local rank
preservation is generally considered implausible. This article relaxes
the assumption of local rank preservation in structural nested models.
The estimators and tests applied in the previous literature are
specific cases of the estimators and tests studied in this article,
with the only difference that some of the estimators in the previous
literature allow for censoring of $Y$. Aside from the issue of
censoring, this article provides a mathematical foundation behind
previous estimators, relaxes the specification of the counterfactual
outcomes as deterministic variables, and allows for a distributional
interpretation of the estimators.

\citet{Enc} conjectures that one can often use standard software to
test whether treatment affects the outcome of interest (without
specifying a model $D_\psi$ for $D$), and to estimate $\psi$.
Lok (\citeyear{Lok,plss}) shows that both testing and estimation can also be
considered from a partial likelihood point of view. As shown in
Lok (\citeyear{Lok, plss}), this approach leads to a subclass of the estimators
and tests studied in this article which can indeed be calculated with
standard software. Example~\ref{Wei} is a specific case of that. The
possibility to use standard software may be a good reason to choose
these estimators in practice. See \citet{Enc} and Lok (\citeyear{Lok,plss}) for
a more elaborate discussion.

The approach adopted in the current article leads to a large class of
estimators and tests. When treatment and covariates change at finitely
many fixed times only, Robins (\citeyear{R93,Rlat}) proposes, without proof, an
optimal procedure for survival and nonsurvival outcomes,
respectively. The optimal choice of estimators or tests under the
framework of this article is another intriguing topic for future
research.

The current article assumes a parametric model $\lambda_\theta$ for
the prediction of treatment changes. In practice, applications have
often used a semi-parametric Cox model for $\lambda_\theta$.
\citet{Lok} shows that specifying $\lambda_\theta$ using a
semiparametric Cox model leads to unbiased estimating equations, which
just as in this article are martingales for the true parameters.
Consistency and asymptotic normality of the resulting estimators for
$D$ still remain to be shown and constitute interesting topics for
future research.

In many applications, observations are censored. \citet{Enc} and
Hern\'{a}n et al. (\citeyear{Her}) have proposed methods to deal with censoring that could
potentially be adapted to the results in this paper. For $D$ of the
form of Example~\ref{DAexa}, \citet{plss} includes proofs with
censoring due to the fact that the study ends, so-called
administrative censoring, using ideas from these previous papers.


The estimators in this paper depend on solving a differential equation
for each observation. In the examples, these equations are simple
enough to be solved analytically. If that is not possible, these
equations should be solved numerically. It might be worth
investigating how a small contamination of the solution to the
differential equation $X(t)$ affects the estimates of treatment
effect.

I conclude with a discussion of the assumptions used in this article.
The most important assumption in this article is the assumption of no
unmeasured confounding (Assumption~\ref{intconf}). As discussed before, this
assumption is valid if all information has been recorded which both
(i) predicts treatment decisions and (ii) is an independent risk
factor for the outcome of interest. The validity of the assumption of
no unmeasured confounding cannot be tested statistically, and depends
on the quality of the recorded information. Therefore, it is for
subject matter experts to decide about the plausibility of the
assumption of no unmeasured confounding. Second, we only estimate the
effect of treatment for which a short duration of treatment has only a
small effect on the distribution of the outcome of interest. The
effect of the treatment on an individual patient may be large, as long
as the probability of such an effect is small for any small duration
of treatment. Third, the assumption of no instantaneous treatment
effect (Assumption~\ref{jump}) is also restrictive: it excludes the estimation of
the effect of treatments that have instantaneous effects, such as
surgery or other point exposures.
The remaining assumptions in this article are mostly benign. The
assumption that the covariate- and treatment process can be
represented by a cadlag process is generally accepted for most medical
situations [see, e.g., Andersen et al. (\citeyear{ABGK})]. The
functions $h_t$ and $D_\psi$ can be chosen such that the regularity
conditions on these functions are satisfied. Even if $h_t$ is not
bounded, it can often be approximated by bounded functions, and results
may follow by a simple application of Lebesgue's dominated convergence
theorem [see, e.g., Example~\ref{Wei}]. The same is true for the
boundedness condition (Assumption~\ref{int}) of the intensity process $\lambda$.
Assumption~\ref{Ycadlag} that the counterfactual process $Y^{(t)}$
is cadlag is impossible to verify, but it is a plausible and
convenient regularity condition.

%
%

\renewcommand{\thetheo}{A.\arabic{theo}}
\setcounter{theo}{0}
\renewcommand{\thecorol}{A.\arabic{corol}}
\setcounter{corol}{2}
\begin{appendix}

\section{Some theory about differential equations}
\label{difeqap}

\begin{theo} \label{isun}
Suppose that a function $D(y,t;\overline{Z}_t)$ satisfies the
following:
\begin{longlist}[(a)]
\item[(a)] (\textup{continuity between the jump times of $Z$}). If
$Z$ does not jump in $(t_1,t_2)$, then
$D(y,t;\overline{Z}_t)$ is continuous in $(y,t)$ on
$[t_1,t_2)$ and can be continuously extended to
$[t_1,t_2]$.
\item[(b)] (\textup{Lipschitz continuity}). For
each $\omega\in\Omega$, there exists a constant $L(\omega)$
such that
\[
|D(y,t;\overline{Z}_t)-D(z,t;\overline{Z}_t)|\leq
L(\omega)|y-z|
\]
for all $t\in[0,\tau]$ and all $y,z$.
\end{longlist}

Suppose, furthermore, that, for each $\omega\in\Omega$, there are no more
than finitely many jump times of $Z$. Then, for each
$t_0\in[0,\tau]$ and $y_0\in\mathbb{R}$, there is a unique
continuous solution $x(t;t_0,y_0)$ to
\[
x'(t)=D(x(t),t;\overline{Z}_t)
\]
with boundary condition $x(t_0)=y_0$ and this solution is
defined on the whole interval $[0,\tau]$.
\end{theo}

This theorem follows from well-known results about
differential equations; see, for example, \citet{anDe}, Chapter~2.

For the next theorem, we also refer to \citet{anDe}, Chapter~2. It is a
consequence of Gronwall's lemma.
\begin{theo} \label{difb} Suppose that $I$ is an open or closed
interval in ${\mathbb R}$,
$f\dvtx I\times {\mathbb R}^n\rightarrow {\mathbb R}^n$ is continuous and
$C\dvtx I\rightarrow[0,\infty)$ is continuous, and suppose that
\begin{equation}
\Vert  f(x,y)-f(x,z)\Vert \leq
C(x) \Vert  y-z\Vert
\label{Lc2}
\end{equation}
for all $x\in I$ and $y,z\in {\mathbb R}^n$.
Then, for every $x_0\in I$ and $y_0\in\mathbb{R}$,
there is a unique solution $y(x)$ of
$y'(x)=f(x,y(x))$ with
$y(x_0)=y_0$, and this solution is defined
for all $x\in I$. If
$g\dvtx I\times {\mathbb R}^n\rightarrow {\mathbb R}^n$ is continuous and
$z\dvtx I\rightarrow {\mathbb R}^n$ is a
solution of
$z'(x)=g(x,z(x))$, then
\begin{eqnarray*}
&&\Vert  y(x)
-z(x)\Vert \\
&&\qquad\leq e^{\int_{x_0}^x C(\xi)d\xi}
\Vert  y(x_0)-z(x_0)\Vert \\
&&\hspace*{8pt}\qquad{}+\int_{x_0}^x e^{\int_{\xi}^x C(\eta)d\eta}
\Vert  f(\xi,z(\xi))-
g(\xi,z(\xi))\Vert \, d\xi
\end{eqnarray*}
for all $x,x_0\in I$ with $x_0\leq x$.
\end{theo}

In \citet{anDe} the interval is always an open interval, but as is
generally known, this can be overcome by extending both $f$ and $g$
outside the closed interval $I$ by taking the values at the boundary
of $I$. This preserves the Lipschitz- and continuity conditions.
Existence and uniqueness on all of finitely many intervals implies
global existence and uniqueness.

We have a differential equation with end condition at
$\tau$, so we are interested in $x,x_0$ with $x\leq x_0$:
\begin{corol}
\label{difbc1}
Suppose that the conditions of Theorem~\textup{\ref{difb}} are satisfied.
Then, for every $x_0\in I$ and $y_0\in\mathbb{R}^n$,
there is a unique solution $y(x)$ of
$y'(x)=f(x,y(x))$ with
$y(x_0)=y_0$, and this solution is defined
for all $x\in I$. If
$g\dvtx I\times {\mathbb R}^n\rightarrow {\mathbb R}^n$ is continuous and
$z\dvtx I\rightarrow {\mathbb R}^n$ is a
solution of
$z'(x)=g(x,z(x))$, then
\begin{eqnarray*}
&&\Vert  y(x)
-z(x)\Vert \\
&&\qquad\leq
e^{\int_{x}^{x_0} C(s)ds}
\Vert  y(x_0)-z(x_0)\Vert \\
&&\qquad\hspace*{7pt}{}+\int_{x}^{x_0} e^{\int_{x}^s C(\eta)d\eta}
\Vert  f(s,z(s))-
g(s,z(s))\Vert \,ds
\end{eqnarray*}
for all $x,x_0$ with $x\leq x_0$.
\end{corol}
\begin{pf*}{Proof} Put
$\tilde{y}(t)=y(x_0-t)$. Then
\[
\tilde{y}'(t)
=-y'(x_0-t)
=-f\bigl(x_0-t,y(x_0-t)\bigr)
=\tilde{f}(t,\tilde{y}(t)),
\]
where $\tilde{f}(t,y)=-f(x_0-t,y)$. Thus,
$\tilde{y}(t)=y(x_0-t)$ is a solution of the
differential equation
$\tilde{y}'(t)=\tilde{f}(t,\tilde{y}(t))$
with boundary condition
$\tilde{y}(0)=y(x_0)=y_0$. Define also
$\tilde{z}(t)=z(x_0-t)$.
Applying Theorem~\ref{difb} on $\tilde{y}$ concludes the proof, as
follows:
\begin{eqnarray*}
\Vert  y(x)-z(x)\Vert
&=& \big\Vert  y\bigl(x_0-(x_0-x)\bigr)-z\bigl(x_0-(x_0-x)\bigr)\big\Vert \\
&=&\Vert  \tilde{y}(x_0-x) - \tilde{z}(x_0-x)\Vert \\
&=&\Vert  \tilde{y}(t)-\tilde{z}(t)\Vert
\end{eqnarray*}
with $t=x_0-x\geq 0.$ Notice that, because of
equation~(\ref{Lc2}),
\[
\Vert  \tilde{f}(t,y)-\tilde{f}(t,z)\Vert \leq
C(x_0-t) \Vert  y-z\Vert =:\tilde{C}(t)
\Vert  y-z \Vert ,
\]
with $\tilde{C}(t)=C(x_0-t)$.
Hence, Theorem~\ref{difb} implies that
\begin{eqnarray*}
\Vert  y(x)-z(x)\Vert
&\leq& e^{\int_0^t\tilde{C}(\xi)d\xi}
\Vert \tilde{y}(0)-\tilde{z}(0)\Vert \\
&&{}+\int_0^t e^{\int_\xi^t \tilde{C}(\eta)d\eta}
\Vert \tilde{f}(\xi,\tilde{z}(\xi))
-\tilde{g}(\xi,\tilde{z}(\xi))\Vert \, d\xi\\
&=&e^{\int_0^t C(x_0-\xi)d\xi}\Vert  y(x_0-0)-z(x_0-0)\Vert \\
&&{}+\int_0^t e^{\int_\xi^t C(x_0-\eta)d\eta}
\Vert \tilde{f}(\xi,\tilde{z}(\xi))
-\tilde{g}(\xi,\tilde{z}(\xi))\Vert \, d\xi.
\end{eqnarray*}
For the first term, we do a change of variables; $\xi$ from $0$ to
$t$, put $s=x_0-\xi$; $d\xi=-ds$. $0\leq\xi\leq t$; $s$ from $x_0-0$
to $x_0-t=x_0-(x_0-x)=x$. We conclude that the first term is equal to
\[
e^{-\int_{x_0}^x C(s)ds}\Vert  y(x_0)-z(x_0)\Vert
= e^{\int_{x}^{x_0} C(s)ds}\Vert  y(x_0)-z(x_0)\Vert .
\]
For the second term, similar changes of variables can be done,
resulting in Corollary~\ref{difbc1}.
\end{pf*}

\section{Mimicking counterfactual outcomes}\label{XTap}

In this appendix we present conditions under which $X(t)$ mimics
$Y^{(t)}$ in the sense that it has the same distribution as $Y^{(t)}$
given $\overline{Z}_t$. This result is used heavily in this article. For the
proofs, which are lengthy and use discretization, we refer
to Lok (\citeyear{Lok,Mimarx}). Section~\ref{XTapT} deals with survival outcomes,
Section~\ref{XTapY} with other outcomes. Survival outcomes require a
different set of assumptions, as will become clear below. The
conditions here are somewhat more restrictive than the ones in
Lok (\citeyear{Lok,Mimarx}), but they are simpler.

\renewcommand{\theassumpt}{B.\arabic{assumpt}}
\setcounter{assumpt}{0}
\renewcommand{\thetheor}{B.\arabic{theor}}
\setcounter{theor}{1}
\renewcommand{\thelemm}{B.\arabic{lemm}}
\setcounter{lemm}{4}
\subsection{Mimicking counterfactual nonsurvival outcomes}
\label{XTapY}

This section contains a sufficient set of regularity conditions to have
existence and uniqueness of a solution $X(t)$ to (\ref{Xdef}),
$
X'(t)=D(X(t),t;\overline{Z}_t)$
with final condition $X(\tau)=Y$, the observed outcome (see
Figure~\ref{Dtfig}). Furthermore, together with Assumption~\ref{inst}
(consistency), they imply that $X(t)$ has the same distribution as
$Y^{(t)}$ given $\overline{Z}_t$.

The regularity conditions below should be read as the following: there exist
conditional distribution functions $F_{Y^{(t+h)}|\overline{Z}_t}$ such
that all these assumptions are satisfied. They can be relaxed to $h$
in a neighborhood of $0$, if this neighborhood does not depend on
$\overline{Z}$. We only consider $h\geq 0$, so the derivative with
respect to $h$ at $h=0$ is always the right-hand derivative.
\begin{assumpt}[(Regularity condition)]\label{sup}
\begin{itemize}
\item (\textit{Support}).
\begin{longlist}[(a)]
\item[(a)] There exist finite numbers $y_1$ and $y_2$
such that all $F_{Y^{(t+h)}|\overline{Z}_t}$ have the same bounded
support $[y_1,y_2]$.
\item[(b)] All $F_{Y^{(t+h)}|\overline{Z}_t}(y)$ have a
continuous nonzero density
$f_{Y^{(t+h)}|\overline{Z}_t}(y)$ on
$y\in[y_1,y_2]$.
\item[(c)] There exists an $\epsilon>0$ such that
$f_{Y^{(t)}|\overline{Z}_t}(y)\geq \epsilon$ for all
$y\in[y_1,y_2]$, $\omega\in\Omega$ and $t\in[0,\tau]$.
\end{longlist}
\item (\textit{Smoothness}). For every $\omega\in\Omega$,
\begin{longlist}[)a)]
\item[(a)] $(y,t,h)\rightarrow F_{Y^{(t+h)}|\overline{Z}_t}(y)$ is
differentiable with respect to $t$, $y$ and $h$ with continuous
derivatives on $[y_1,y_2]\times[t_1,t_2)\times
{\mathbb R}$ if $Z$ does not jump in $(t_1,t_2)$, with a
continuous extension to $[y_1,y_2]\times[t_1,t_2]\times
[0,\infty)$.
\item[(b)] The derivatives of $F_{Y^{(t+h)}|\overline{Z}_t}(y)$
with respect to $y$ and $h$ are bounded by constants $C_1$ and $C_2$,
respectively.
\item[(c)] $\frac{\partial}{\partial y} F_{Y^{(t)}|\overline{Z}_t}(y)$ and
$\frac{\partial}{\partial h}|_{h=0}
F_{Y^{(t+h)}|\overline{Z}_t}(y)$ have derivatives with respect to $y$
which are bounded by constants $L_1$ and $L_2$, respectively.
\end{longlist}
\end{itemize}
\end{assumpt}

The support conditions may be restrictive for certain
applications. Nevertheless, most real-life situations can be
approximated this way, since $y_1$ and $y_2$ are unrestricted and
$\epsilon>0$ is unrestricted, too. Although the support conditions may well
be stronger than necessary, they simplify the analysis considerably
and, for that reason, they are adopted here.
The smoothness conditions allow for nonsmoothness where the
covariate- and treatment process $Z$ jumps. This is important, since if
the covariate- and treatment process $Z$ jumps, this can lead to a
different prognosis for the patient and thus to nonsmoothness of the
functions concerned.
\begin{theor}[\textup{(Mimicking counterfactual outcomes)}] \label{thm1}
Suppose that regularity Condition~\textup{\ref{sup}} is satisfied. Then
$D(y,t;\overline{Z}_t)$ exists. Furthermore, for every
$\omega\in\Omega$, there exists exactly one continuous solution $X(t)$
to $X'(t)=D(X(t),t;\overline{Z}_t)$ with final condition
$X(\tau)=Y$. If also Assumption~\textup{\ref{inst}} (consistency) is
satisfied and there are no more than finitely many times $t$ for which
the probability that the covariate- and treatment process jumps at
$t$ is greater than $0$, then this $X(t)$ has the same distribution as
$Y^{(t)}$ given $\overline{Z}_t$ for all $t\in[0,\tau]$.
\end{theor}

For a proof we refer to Lok (\citeyear{Lok,Mimarx}).

\subsection{Mimicking counterfactual survival outcomes}
\label{XTapT}

This section contains a sufficient set of regularity conditions to have
existence and uniqueness of a solution $X(t)$ to equation~(\ref{Xdef}),
$X'(t)=D(X(t),t;\overline{Z}_t)$
with final condition $X(\tau)=Y$, the observed outcome (see
Figure~\ref{Dtfig}). Furthermore, together with Assumption~\ref{inst}
(consistency) and Assumptions~\ref{ttcons} and~\ref{ttinst} below,
they imply that $X(t)$ has the same distribution as $Y^{(t)}$ given
$\overline{Z}_t$. The conditions here are natural conditions if the
outcome of interest $Y$ is a survival time.

\renewcommand{\theassumpt}{B.\arabic{assumpt}}
\setcounter{assumpt}{2}

As compared to Section~\ref{XTapY}, we make two extra assumptions. The
first is a consistency assumption, stating that stopping treatment
after death does not change the survival time. The second assumption
states that there is no instantaneous effect of treatment at the time
the patient died [notice that the difference between
$Y^{(Y)}$, the outcome with treatment stopped at the
survival time $Y$, and $Y$ is in treatment \textit{at} $Y$].
\begin{assumpt}[(Consistency)] \label{ttcons}
$Y^{(t)}=Y$ on $\{\omega\dvtx Y< t\}\cup\{\omega\dvtx Y^{(t)}<t\}$.
\end{assumpt}
\begin{assumpt}[(No instantaneous effect of treatment at the
time the patient died)] \label{ttinst} $Y^{(t)}=Y\mbox{ on }\{\omega\dvtx Y=
t\}\cup\{\omega\dvtx Y^{(t)}=t\}.$
\end{assumpt}

Under these assumptions, treatment in the future does not cause
or prevent death at present or before:
\begin{lemm} \label{conslem} Under Assumptions~\textup{\ref{ttcons}}
and~\textup{\ref{ttinst}}:
\begin{longlist}[(a)]
\item[(a)] For all $h\geq 0$: $Y^{(t+h)}=Y$ on $\{\omega\dvtx Y\leq t\}\cup\bigcup_{h\geq
0}\{\omega\dvtx Y^{(t+h)}\leq t\}$.
\item[(b)] For all $(y,t,h)$ with $y\leq t+h$ and $h\geq 0\dvtx \{\omega\dvtx Y^{(t+h)}\leq y\}=\{\omega: Y\leq
y\}$.
\end{longlist}
\end{lemm}

For a proof we refer to \citet{Lok}.

\renewcommand{\theassumpt}{B.\arabic{assumpt}}
\setcounter{assumpt}{5}

If the outcome is survival, the support condition in Assumption~\ref{XTapY},
saying that all $F_{Y^{(t+h)}|\overline{Z}_t}$ have the same bounded
support $[y_1,y_2]$, will not hold. The reason for this is
as follows. $\overline{Z}_t$ includes the covariate-measurements and
treatment until time $t$. If covariates and treatment were measured at
time $t$, it cannot be avoided to include in $\overline{Z}_t$ whether
or not a patient was alive at time $t$. Given that a patient is dead at
time~$t$ and given his or her survival time, the distribution of this
survival time cannot have the fixed support $[y_1,y_2]$,
which is independent of $t$. Also, given that a patient is alive at time
$t$, this is hardly ever the case; one often expects that $t$ is the left limit
of the support.  Thus, in case the outcome is survival, the support
condition for Theorem~\ref{thm1} has to be slightly changed.
\begin{assumpt}[(\textit{Support})]\label{tsup}
 There exists a finite number $y_2\geq\tau$ such that:
\begin{longlist}[(a)]
\item[(a)] For every $\omega\in\Omega$ and $t$ with $Y>t$, all
$F_{Y^{(t+h)}|\overline{Z}_t}$ for $h\geq 0$ have support
$[t,y_2]$.
\item[(b)] For every $\omega\in\Omega$ and $t$ with $Y>t$, all
$F_{Y^{(t+h)}|\overline{Z}_t}(y)$ for $h\geq 0$ have a
continuous nonzero density
$f_{Y^{(t+h)}|\overline{Z}_t}(y)$ on
$y\in[t+h,y_2]$.
\item[(c)] There exists a number $\epsilon>0$ such that, for all
$\omega\in\Omega$ and $t$ with $Y>t$,
$f_{Y^{(t)}|\overline{Z}_t}(y)>\epsilon$ for
$y\in[t,y_2]$.
\end{longlist}
\end{assumpt}

Next we look at the differentiability conditions in Assumption~\ref{sup}. It does
not seem reasonable to assume that
$F_{Y^{(t+h)}|\overline{Z}_t}(y)$ is continuously
differentiable with respect to $h$ and $y$ on
$(h,y)\in[0,\infty)\times [t,y_2]$
since, for $y\leq t+h$,
$F_{Y^{(t+h)}|\overline{Z}_t}(y)=F_{Y|\overline{Z}_t}(y)$
[Lemma~\ref{conslem}(b)]. Therefore, the derivative of
$F_{Y^{(t+h)}|\overline{Z}_t}(y)$ with respect to $h$ is
likely not to exist at $y=t+h$ (and is equal to zero for
$y<t+h$). Also, the derivative of
$F_{Y^{(t+h)}|\overline{Z}_t}(y)$ with respect to $y$ may
not exist at $y=t+h$, because of the different treatment before and
after $t+h$. For survival outcomes, we replace the smoothness conditions
of Assumption \ref{sup} by the following:
\begin{assumpt}[(\textit{Smoothness})]\label{tsmooth}
For every $\omega\in\Omega$:
\begin{longlist}[(a)]
\item[(a)] If $Z$ does not jump in $(t_1,t_2)$ and $Y>t_1$, the
restriction of $(y,t,h)\rightarrow
F_{Y^{(t+h)}|\overline{Z}_t}(y)$ to
$\{(y,t,h)\in[t_1,y_2]\times[t_1,t_2)\times
{\mathbb R}_{\geq 0}\dvtx y\geq t+h\}$ is $C^1$ in
$(y,t,h)$.
\item[(b)] The derivatives of $F_{Y^{(t+h)}|\overline{Z}_t}(y)$
with respect to $y$ and $h$ are bounded by constants $C_1$ and $C_2$,
respectively, for $y\in[t+h,y_2]$.
\item[(c)] $\frac{\partial}{\partial y} F_{Y^{(t)}|\overline{Z}_t}(y)$ and
$\frac{\partial}{\partial h}|_{h=0}
F_{Y^{(t+h)}|\overline{Z}_t}(y)$ have derivatives with respect to $y$
which are bounded by constants $L_1$ and $L_2$, respectively,
for $y\in[t+h,y_2]$.
\end{longlist}
\end{assumpt}
\renewcommand{\thetheor}{B.\arabic{theor}}
\setcounter{theor}{8}

The smoothness condition above concentrates on $y\geq t+h$. For
$y\in[t,t+h)$ we can choose
$F_{Y^{(t+h)}|\overline{Z}_t}(y)=F_{Y|\overline{Z}_t}(y)$ because of
Lemma~\ref{conslem}(b). Because of Assumption~\ref{inst} (consistency),
$F_{Y|\overline{Z}_t}$ has the same support as
$F_{Y^{(\tau)}|\overline{Z}_t}$, so $F_{Y|\overline{Z}_t}$ has support
$[t,y_2]$ if $Y>t$ [Assumption~\ref{tsup}(a)]. Assume the following
\begin{assumpt}[(\textit{Smoothness})]\label{Ytsup}
For all $\omega\in\Omega$ and $t$ with $Y>t$,
$F_{Y|\overline{Z}_t}(y)$ is continuous and strictly increasing on its
support $[t,y_2]$.
\end{assumpt}
\begin{theor}[(\textup{Mimicking counterfactual survival outcomes})] \label{tthm1}
Suppose that regularity Conditions \textup{\ref{tsup}, \ref{tsmooth}}
and~\textup{\ref{Ytsup}} are satisfied. Then $D(y,t;\overline{Z}_t)$
exists. Furthermore, for every $\omega\in\Omega$, there exists exactly
one continuous solution $X(t)$ to
$X'(t)=D(X(t),t;\overline{Z}_t)$ with final condition
$X(\tau)=Y$. If also Assumptions~\textup{\ref{inst}, \ref{ttcons}}
and~\textup{\ref{ttinst}} (consistency and no instantaneous treatment effect at
time of death) are satisfied, then this $X(t)$ has the same
distribution as $Y^{(t)}$ given $\overline{Z}_t$ for all
$t\in[0,\tau]$.
\end{theor}

For a proof we refer to Lok (\citeyear{Lok,Mimarx}).

\renewcommand{\thelemma}{C.\arabic{lemma}}
\setcounter{lemma}{0}

\section{Measurability issues}
\label{smtb}

In most of this article we assume that the function which maps
$(X(t),\overline{Z}_{t-})$ to $X(t_0)$, with $t_0<t$, is a
measurable function on $\mathbb{R}\times \overline{\mathcal{Z}}_{t-}$,
with the projection $\sigma$-algebra on $\overline{\mathcal{Z}}_{t-}$
(see Section~\ref{sets}). Moreover, we sometimes assume that the
function which maps $(X(t_0),\overline{Z}_{t-})$ to $X(t)$,
with $t_0<t$, is a measurable function on $\mathbb{R}\times
\overline{\mathcal{Z}}_{t-}$. In this appendix we give sufficient
conditions for this. If these two functions are measurable,
$\sigma(\overline{Z}_t,X(t))$ is a filtration, and, moreover,
$\sigma(\overline{Z}_t,X(t))$ is the same as
$\sigma(\overline{Z}_t,X(0))$
[see equation~(\ref{fil}) in Section~\ref{mgs}].
\begin{lemma} \label{mtb} Suppose that $D$ satisfies regularity
  Assumption~\textup{\ref{Dreg}} and that, for each $\omega\in\Omega$, $Z$ jumps
  at most finitely many times. Then the function which maps
  $(X(t),\overline{Z}_{t-})$ to $X(t_0)$, with $t_0<t$, is
  a measurable function from $\mathbb{R}\times \overline{\mathcal{Z}}_{t-}$ to $\mathbb{R}$.
  Also, the function which maps
  $(X(t_0),\overline{Z}_{t-})$ to $X(t)$, with $t_0<t$, is
  a measurable function from $\mathbb{R}\times \overline{\mathcal{Z}}_{t-}$ to $\mathbb{R}$.
\end{lemma}

For the proof of this result, which is quite technical since many
results on differential equations are nonconstructive, we refer to
\citet{Lok}. The proof uses the idea behind Euler's forward method to
approximate the solution to the differential equation.

\renewcommand{\thelemmas}{D.\arabic{lemmas}}
\setcounter{lemmas}{0}
\section{A convergence result}
\label{GCap}

The following lemma is a worked-out case of theory from \citet{Vaart},
Chapter~19.
\begin{lemmas}\label{pGC} Let $X_1,X_2,\ldots$ be i.i.d. random
variables with values in a measurable space $\mathcal{X}$.  Let
$\{f_{\theta}\dvtx\theta\in\Theta\}$ be a collection of
measurable functions from $\mathcal{X}$ to $\mathbb{R}^k$ indexed by a
subset $\Theta\subset \mathbb{R}^d$ which contains an open
neighborhood $\Theta_0$ of $\theta_0$. Suppose that
$\theta\rightarrow f_\theta(x)$ is continuous on $\Theta_0$
for every $x\in \mathcal{X}$. Suppose also that there exists a measurable
function $F$ on $\mathcal{X}$ such that $\| f_\theta\|\leq F$
for every $\theta\in\Theta_0$ and such that $E F(X_1)$
exists.
Then if $\hat{\theta}$ converges in probability to $\theta_0$,
\[
P_n f_{\hat{\theta}} \mathop{\stackrel{P}{\rightarrow}} Ef_{\theta_0}(X_1),
\]
where $P_n$ indicates the empirical distribution of
$X_1,X_2,\ldots,X_n$.
\end{lemmas}
\begin{pf*}{Proof} Notice that
\[
\| P_nf_{\hat{\theta}} - Ef_{\theta_0}(X_1)\|
\leq \| P_nf_{\hat{\theta}} - Ef_{\hat{\theta}}(X_1)\|
+\| Ef_{\hat{\theta}}(X_1)-
Ef_{\theta_0}(X_1)\|.
\]
We show that both terms converge to zero in probability.
Choose $\Theta_1\subset\Theta_0$ compact and such that it contains an
open neighborhood of $\theta_0$. Example~19.8 from \citet{Vaart}
implies that, under the conditions above,
\[
\sup_{\theta\in\Theta_1}
\| P_n f_\theta -Ef_\theta(X_1)\|
\rightarrow 0\qquad \mbox{a.s.}
\]
Since $\hat{\theta}\mathop{\stackrel{P}{\rightarrow}}\theta_0$ and $\Theta_1$ contain an
open neighborhood of $\theta_0$, this implies that the first term
converges in probability to zero.  For the second term, notice that, on
$\Theta_0$, $\theta\rightarrow f_{\theta}(x)$ is continuous in
$\theta$ and that each of the components of $f_\theta$ is bounded by
the integrable function $F$,
so that Lebesgue's dominated convergence theorem implies that
$Ef_\theta(X_1)$ is continuous in $\theta$ on
$\Theta_0$. Thus, since $\hat{\theta}\mathop{\stackrel{P}{\rightarrow}}\theta_0$, also the
second term converges in probability to zero.
\end{pf*}
\end{appendix}

\section*{Acknowledgments}
An earlier version of this
  article was a chapter of my Ph.D.~dissertation at the Free
  University of Amsterdam. I am indebted to Richard Gill and Aad van
  der Vaart for their support, insight and encouragement on this
  project. I also thank James Robins for fruitful discussions. I thank
  both referees for helpful suggestions. Besides that, I thank the
  Netherlands Organization for Scientific Research (NWO) for a
  scholarship. 

\printaddresses

\end{document}